\documentclass[final,12pt]{elsarticle}

\usepackage[margin=4.46cm]{geometry}
\usepackage{graphicx}
\usepackage{tikz}
\usepgflibrary{shapes}
\usepackage{epsfig}
\usepackage{amsmath}
\usepackage{amssymb}
\usepackage{enumitem}

\usepackage{caption}
\usepackage{colortbl} 

\def \Rset {{\mathbb R}}

\def \Zset {{\mathbb Z}}

\newcommand{\nit}{\noindent}
\newcommand{\no}{\nonumber}
\newcommand{\be}{\begin{equation}}
\newcommand{\ee}{\end{equation}}
\newcommand{\ba}{\begin{eqnarray}}
\newcommand{\ea}{\end{eqnarray}}
\newcommand{\bi}{\begin{itemize}}
\newcommand{\ei}{\end{itemize}}
\newcommand{\br}{\begin{eqnarray}}
\newcommand{\er}{\end{eqnarray}}

\newcommand{\real}{\mathbb{R}}

\RequirePackage[latin1]{inputenc}
\RequirePackage{hyperref}
\RequirePackage{amsmath,amssymb,amsthm, amsfonts,mathtools}
\theoremstyle{plain}

\newtheorem{thm}{Theorem}[section]

\newtheorem{dfn}[thm]{Definition}
\newtheorem{lem}[thm]{Lemma}

\newtheorem{rmk}[thm]{Remark}
\newtheorem{ques}{Question}
\usepackage{xcolor}
\usepackage{amsmath}
\usepackage{amssymb}
\usepackage{mathtools}
\usepackage{amsthm}
\usepackage{hyperref}
\usepackage{verbatim}
\usepackage{svg}

\usepackage{amsmath,amsfonts,bm}

\def\eqref#1{equation~(\ref{#1})}

\def\1{\bm{1}}

\DeclareMathAlphabet{\mathsfit}{\encodingdefault}{\sfdefault}{m}{sl}
\SetMathAlphabet{\mathsfit}{bold}{\encodingdefault}{\sfdefault}{bx}{n}

\usepackage{taost}
\usepackage{listings}

\journal{\hspace{1.0 in} \textup{1}}

\begin{document}

\begin{frontmatter}



\title{Lagrangian, Game Theoretic and PDE Methods for Averaging 
G-equations in Turbulent Combustion: Existence and Beyond}


\author{Jack Xin$^{a}$, Yifeng Yu$^{a}$, Paul Ronney\,$^{b}$}
                     
\address[J. Xin, Y. Yu]{Department of Mathematics, University of California, {Irvine},CA 92697, {USA}}
            
\address[P. Ronney]{Department of Aerospace and Mechanical Engineering, University of Southern California, LA, {CA 90089}, {USA}}
            
\begin{abstract}
G-equations are popular level set Hamilton-Jacobi nonlinear partial differential equations (PDEs) of first or second order arising in turbulent combustion. 
Characterizing the effective burning velocity (also known as the turbulent burning velocity) is a fundamental problem there.
We review relevant studies of the G-equation models with a focus on
both the existence of effective burning velocity (homogenization), and 
its dependence on physical and geometric parameters (flow intensity and curvature effect) through representative examples. The corresponding physical background is also presented to  provide motivations for mathematical problems of interest.

The {\it lack of coercivity} of  Hamiltonian is a hallmark of G-equations.  
 When either the curvature of the level set or  the strain effect of fluid flows is accounted for, 
the Hamiltonian  
becomes {\it highly non-convex and nonlinear}.
In the absence of 
coercivity and convexity, PDE (Eulerian) approach 
suffers from insufficient compactness to establish averaging (homogenization). 
We review and illustrate a suite of Lagrangian tools, 
most notably min-max (max-min) game representations of curvature and strain G-equations, working in tandem with analysis of streamline structures of fluid flows and PDEs. 
We discuss open problems for future development in this emerging area of dynamic game 
analysis for averaging non-coercive, non-convex, and nonlinear PDEs such as geometric (curvature-dependent) PDEs with advection. 
\medskip

{\it To appear in Bulletin of the American Math Society, 2024.}
\end{abstract}
 
\thispagestyle{empty}

\end{frontmatter}


\setcounter{page}{1}
\date{}
\section{Introduction}
\label{intro}
Turbulent combustion is employed in practically all mobile and stationary power generation devices because turbulence increases the mass consumption rate of reactants to values much greater than laminar flames can achieve.  This in turn increases the heat release rate and thus power available from a gas turbine combustor, reciprocating piston engine, or rocket motor of a given size.  Few combustion engines would function without the increase in mixing and burning rates caused by turbulence.  A fundamental and  most practically important property of a turbulent flame in a pre-mixture of reactants is how the effective propagation speed ($s_T$, also called turbulent flame speed) is affected by turbulence intensity among other factors including local flame front curvature, hydrodynamic strain, combustor geometry, thermal expansion, heat losses, and turbulence energy spectrum \cite{Brad_92,Pet_00,Dris_08}.  Most of our understanding of these effects to date in combustion science has been developed from the so-called ``flamelet'' models that presume the existence of a continuous flame sheet that is wrinkled by turbulence and thereby affecting the flame area which in turn affects the mass consumption rates. The local consumption rate per unit area corresponds to the propagation rate of a planar front in the same mixture.
\medskip
\begin{figure}
\begin{center}
\includegraphics[scale=0.8]{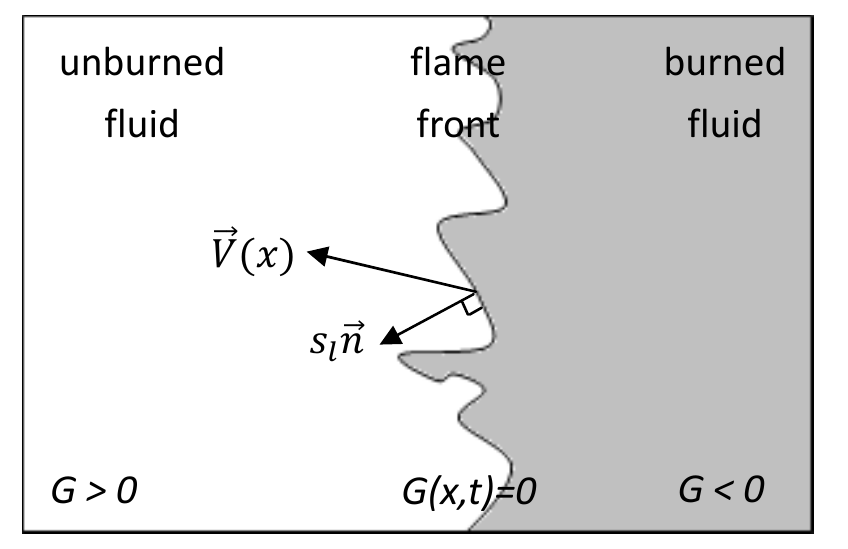}
\caption{ Level set modeling of flame propagation.}
\label{G-modeling}
\end{center}
\end{figure}

A well-known flamelet type model in turbulent combustion is the so called G-equation  \cite{Mark_64,Will_85,Will_85b,Pet_00}, which takes the following form:
\be\label{main-eq}
G_t+s_l\, |DG|+V(x)\cdot DG=0.
\ee
Here $G=G(x,t)\in C(\Rset^n\times (0,\infty), \Rset)$ is a reference function whose zero level set $\{x\in \Rset^n|\ G(x,t)=0\}$ represents flame front at time $t$; $\{x\in \Rset^n|\ G(x,t)>0\}$ is the  unburned area and $\{x\in \Rset^n|\ G(x,t)<0\}$ the burned area;  $V\in C(\Rset^n, \Rset^n)$  is the velocity field of the ambient fluid (e.g. blowing wind in a wild fire or stirring a mixture of air and gasoline inside car engines),  see Fig. \ref{G-modeling} below; $D$ the spatial gradient. The motion of the flame front is driven by chemical reaction represented by the so called laminar flame speed $s_l$ (or local burning velocity) and the ambient fluid. Hence the flame front obeys the motion law: 
\be
{v}_{\vec{n}}=s_l+V(x)\cdot \vec{n}, \label{mlaw}
\ee
i.e., the front speed ${v}_{\vec{n}}$ in the normal direction $\vec{n}$ is the laminar flame speed ($s_l$) plus the projection of fluid velocity $V$ along the normal.
The $s_l$ might not be a constant in general, and its geometric or physical modeling leads to different types of 
G-equations.  
\medskip

The ratio of flame thickness and Kolmogorov scale (the small scale in fluid flows below which dissipation dominates) determines different combustion regions. G-equation is most suited in the so called corrugated flamelet regime where the ratio is less than one \cite{Pet_00}. 
Below we derive G-equation, assuming that the $G$ function and the flame front are smooth. Fix $t>0$ and  a point $x$ on the flame front $\Gamma_t$.  Let $x(s)$ be the flame propagation route starting from $x(t)=x$ for $s\geq t$. Then $G(x(s),s)=0$ and  the motion law (\ref{mlaw}) implies that
\be\label{normal-speed}
\dot x(s)={v}_{\vec{n}(s)}\, \vec{n}(s)=
(s_l +V(x(s))\cdot \vec{n}(s))\, \vec{n}(s),
\ee
where $\vec{n}(s)$ is the outward normal vector of $\Gamma_s$ at $x(s)$ that can be expressed as $\vec{n}(s)=DG(x(s),s)/ |DG(x(s),s)|$. By chain rule: 
$$
0={dG(x(s),s)\over ds}=G_t+DG\cdot \dot x(s)=0.
$$
Plugging in (\ref{normal-speed}) immediately gives G-equation (\ref{main-eq}), which was formally introduced by Williams \cite{Will_85,Will_85b} in 1985 although the above derivation and  an earlier form of G-equation first appeared in Markstein's work 
\cite{Mark_64} in 1964.  Also,  we note that the  G-equation model and crystal growth front propagation model served as two major  sources in the systematic mathematical development of the level-set method by Osher and Sethian \cite{OS_88}.  

Two observations are in order.
\medskip

(1) The reference function $G(x,t)$ has no physical meaning except its zero level set (called nonreactive in combustion literature).  Solutions to G-equations are in general not $C^1$ although the equation is derived under smoothness assumptions. Mathematically, the solutions need to be defined in the sense of viscosity solutions whose basic definition and backgrounds are in section 2 providing mathematical foundation of level set equations.  In particular, the choice of different reference functions does not impact the zero level set. Precisely speaking,  if both $\tilde G$ and $G$ are solutions to equation (\ref{main-eq}) subject to 
$$
\{x\in \Rset^n|\ \tilde G(x,0)\}=\{x\in \Rset^n|\ G(x,0)\}
$$
then 
$$
\{x\in \Rset^n|\ \tilde G(x,t)=0\}=\{x\in \Rset^n|\ G(x,t)=0\}
$$
for all $t>0$.  An intuitive  way to understand this is that $G$ is a solution if and only if  $\tilde G=f(G)$  is also a solution for any differentiable function $f:\Rset\to \Rset$ satisfying $f'>0$. See \cite{CGG1991,EVS1991} for rigorous proofs. 

\medskip

(2)  In more comprehensive combustion models, fluid velocity and chemical transport are determined by a coupled system of Navier-Stokes and
reaction-diffusion-advection equations, which  is called the direct method. A flamelet approach based on G-equations neglects the effects of thermal expansion (and so the above two-way coupling between the flame and turbulent flow field) as well as heat losses but may retain curvature and strain dependence on $s_l$ in the G-equation formulation.  In the simplified setting of G-equation, the fluid velocity (turbulent flow field) is prescribed with the focus to understand how the flame propagates in an ambient fluid. Such an approach is known as passive modeling, making it more feasible 
to theoretically study (so called ``pencil and paper") or compute the impacts of significant factors in combustion, e.g.,  flow perturbation, flame intensity and flame stretch. That the passive modeling  
may yield physically meaningful results has been 
demonstrated by experiments \cite{Ronney_95,Shy_96,Abid_99} that employ aqueous autocatalytic chemical reactions rather than gaseous flames because the thermal expansion and heat loss effects are absent in the aqueous mixtures.

To initiate a practical simulation of G-equation for turbulent flows, solutions of Navier-Stokes equations are often computed first to set up the velocity field $V$. For numerical approximations of G-equations and applications in combustion, we refer to \cite{ZR_94,  OS_88,OF_02,LXY_13a,KLX_21, LLM2021, EnsKin_23} among others.

\medskip

\subsection{Different types of G-equations} In general, the local burning velocity $s_l$ depends on flame stretch \cite{Pet_00}: 
\be\label{eq:generalspeed}
s_L=(s_{L}^{0}-d\,s_{L}^{0}\, \kappa+d\, \vec {n}\cdot S \cdot  \vec {n})_{+}.
\ee
See section 2.6 in \cite{Pet_00} for a formal derivation of relevant expressions from the reaction-diffusion-advection equations. Here

$\bullet$ $s_{L}^{0}$ is a positive constant representing the laminar flame speed of the un-stretched flame. The positive  part $(a)_{+}=\max\{a,0\}$ is imposed to avoid negative laminar flame speed since materials 
cannot be ``unburned". This correction usually is not explicitly mentioned in combustion literature since,  by default, the laminar flame speed is always assumed to be non-negative there. However mathematically, large positive curvature $\kappa$ or negative strain rate $\vec {n}\cdot S \cdot  \vec {n}$ could occur as time evolves. Hence it is necessary to explicitly add $(\cdot)_{+}$ if physical validity is taken into consideration in the modeling of flame propagation \cite{ZR_94,plus2_97};

$\bullet$ the expression $-d\,s_{L}^{0} \kappa+d\, \vec {n}\cdot S \cdot  \vec {n}$ represents the correction due to flame stretching. The constant $d$ is the Markstein length which is proportional to laminar flame thickness. Intuitively speaking,  the local burning velocity is affected by chemical reactions (or burning temperature) and heat release that are related to the surface area of the flame front. 
Suppose that a smooth hypersurface in $\Rset^n$ is moving in the velocity field $V$. Then the front surface stretch rate is
$$
{1\over \sigma}{d\sigma\over dt}=\mathrm{div}(V)-\vec {n}\cdot S \cdot  \vec {n}.
$$
Here $\sigma$ is the surface element area. The derivation of the above formula could be found in \cite{Mat_83, LXY_13a}.

   The first correction term $-\, d\, s_{L}^{0}\, \kappa$ in (\ref{eq:generalspeed}) describes the curvature effect \cite{M1951, Mark_64, Pet_00}, where   
$$
\kappa=\mathrm{div} \left({DG\over |DG|}\right)
$$
is the mean curvature of the flame surface. 

The second correction  term $d\, \vec {n}\,\cdot S\, \cdot \, \vec {n}$ in  (\ref{eq:generalspeed}) quantifies the straining effect on the flame from the flow \cite{Mat_83,LXY_13a}. The normal vector $\vec{n}$ at the flame surface points in the direction of the unburned region, and
$$
S={DV+(DV)^{\top}\over 2}
$$
is called the strain rate tensor.

For convenience, we set the constant 
$$s_{L}^{0}=1
$$
throughout this paper. 

If (\ref{eq:generalspeed}) is plugged into the equation (\ref{main-eq}), the general G-equation is a 2nd order quasilinear parabolic equation:
\be\label{eq:generalG}
G_t+H(D^2G, DG, x)=0,
\ee
where the associated Hamiltonian is
$$
H(M,p,x)=\left(|p|-d\left(\mathrm{tr}M-{p\cdot M\cdot p\over |p|^2}-{p\cdot S \cdot p\over |p|}\right)\right)_{+}+V(x)\cdot p
$$
for $(M,p)\in S^{n\times n}\times \left(\Rset^n\backslash\{0\}\right)$. At $p=0$, following \cite{CIL_92}, $H(M,0,x)$ is set to $0$ for subsolutions and $2dn||A||$ for supersolutions. Here $S^{n\times n}$ is the space of $n\times n$ real symmetric matrices.  Note that $H$ is monotonically non-increasing with respect to $M$, non-coercive ($H\not \rightarrow +\infty$ as $|p| \rightarrow +\infty$ if 
$|V|>1$) and non-convex with respect to the $p$ variable if $d>0$.
\medskip 

For simplicity, we will treat curvature and strain terms separately. 

\subsubsection{Inviscid G-equation ($d=0$)} In the basic setting,  the local burning velocity $s_l$ is taken to be  a constant ($s_l=s_L^0=1$), i.e., ignore the curvature effect and  strain rate by letting $d=0$.  The corresponding G-equation is then a first-order convex but in general non-coercive (if $|V|>1$) Hamilton-Jacobi equation (HJE):
$$
G_t+|DG|+V(x)\cdot DG=0.
$$
In this situation,  the laminar flame speed could depend on locations, i.e., $s_l=a(x)$ for some positive function $a(x)$. A notable example is the celebrated Rothermel surface fire spread model for forest fire \cite{Rothermal}, where 
$a(x)$ depends on
the landscape (uphill or downhill) and the types of plants there. Due to diffusion in a numerical scheme or artificial viscosity to stabilize computation \cite{KAW1988,FIL_97,OF_02}, it is of interest to assess such effects by considering the viscous version of G-equation:
\be
-d\, \Delta G+G_t+|DG|+V(x)\cdot DG=0,
\label{vgeq}
\ee
where the Laplacian term mimics numerical diffusion. 
The Laplacian term is also found to be relevant for 
representing the effect of thermal relaxation under
transverse heat diffusion in the preheat zone of a wrinkled
front \cite{Clavin_85,FIL_97}. Moreover, the viscous G-equation (\ref{vgeq}) is adopted for 
dynamic modeling of large eddy simulation of turbulent premixed combustion in homogeneous isotropic turbulence in \cite{FIL_97}.
Surprisingly, such a second order regularization at any small $d > 0$ acts as a singular perturbation drastically altering the effective burning velocity, see \cite{LXY_11} for a mathematical proof and \cite{FIL_97} for numerical observations based on Fig. 1 therein.

\subsubsection{Curvature G-equation}  Geometry on the flame front provides a phenomenological way to describe the change of temperature along the flame front, that impacts the rate of chemical reactions and hence the local burning speed.  Intuitively, it says that if the flame front bends 
toward the cold region (un-burned area, point C in Fig. \ref{curv} below), 
the flame propagation slows down. If the flame front bends toward the hot 
spot (burned area, point B in Fig. \ref{curv}), it burns faster.

\begin{figure}
\begin{center}
\includegraphics[scale=0.8]{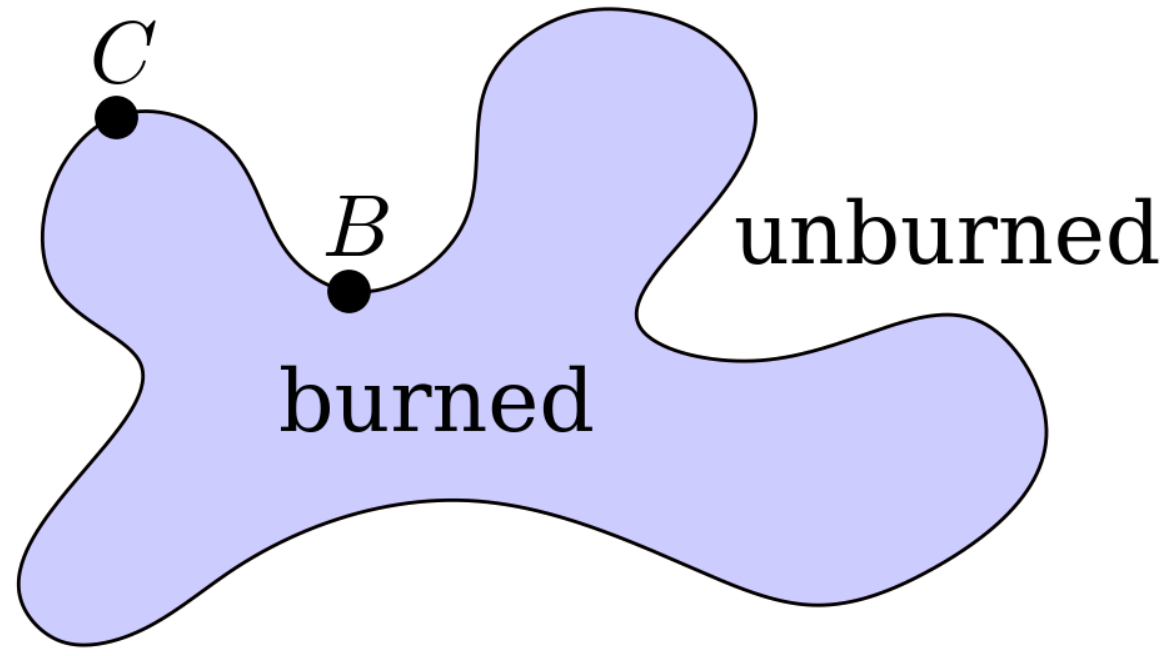}
\caption{Curvature effect: faster burn at B (bending towards hot region) than C.}
\label{curv}
\end{center}
\end{figure}

There is a vast literature in turbulent combustion literature discussing the impact of curvature effect. Below is the most recognized empirical linear relation first proposed by Markstein \cite{M1951, Mark_64} (see also \cite{Pet_00}):
\be\label{laminar}
s_l=s_{L}^{0}(1- d\; \kappa)_{+}.
\ee
Here the strain rate term in (\ref{eq:generalspeed}) is ignored. This correction leads to a mean curvature type  equation with advection:
\be
G_t + \left(1- d\, \, \mathrm{div}\left({DG\over |DG|}\right)\right)_{+}|DG|+V(x)\cdot DG=0.\label{ge1}
\ee
This equation has both non-convex and non-coercive nature. Moreover, due to the presence of the physical cut-off $()_+$, this quasilinear parabolic equation is highly  degenerate, causing non-existence of effective burning velocity in 3D flows (see a later section).  

\medskip

\subsubsection{Strain G-equation} By setting $s_{L}^{0}=1$, the expression for laminar speed is 
\be
s_L=(1+d\, \vec {n}\cdot S \cdot  \vec {n})_{+}. \label{claw}
\ee
This leads to a nonconvex-noncoercive first order HJE: 
\be\label{strainG}
G_t+\left(1+d\ {DG\cdot S(x)\cdot DG\over |DG|}\right)_+|DG|+V(x)\cdot DG=0.
\ee
Nonconvex first order Hamilton-Jacobi equations first appeared in zero sum two-person differential games \cite{I1965,ES1984}. The strain G-equation provides a natural physical example. Although there is a large volume of combustion literature investigating the strain rate effect, there is not much mathematical work studying the above PDE, a nonconvex HJE with physical meaning. 

\subsection{Turbulent burning velocity and homogenization of G-equation}

\begin{figure}
\begin{center}
\includegraphics[scale=0.8]{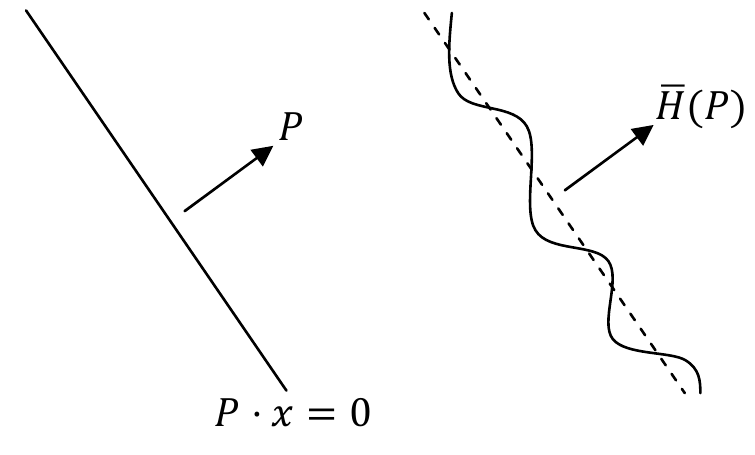}
\caption{A classical (left) vs. an effective (right, dashed) planar front by averaging an oscillatory flame front (right, solid).}
\label{average-plane}
\end{center}
\end{figure}

Effective burning velocity (a.k.a turbulent flame speed)
is one of the most fundamental quantities in turbulent combustion.  Roughly speaking, it is the flame propagation speed after averaging the fluctuations around the mean flame front. Establishing its existence rigorously has in general remained an open problem.  Below are two equivalent perspectives to describe it from the 
G-equation.
   Fix a unit direction $p\in \Rset^n$. Assume that the initial flame front is
   planar $\{x\in \Rset^n|\ G(x,0)=x\}=\{x\in \Rset^n|\ p\cdot x=0\}$. Due to the motion of the ambient fluid, the flame front wrinkles as time evolves, see 
   Fig. \ref{average-plane}. 

   \medskip
   
   {\bf $\bullet$  Perspective 1:} The mean flame front is still planar and might eventually propagate with a steady speed $s_T(p)$.

Mathematically, this implies that 
\be\label{steady-form}
G(x,t)\approx p\cdot x+v(x)-t\, s_T(p),
\ee
where the right hand side is an approximate  traveling wave solution of the G-equation.  The $s_T(p)$, if it exists, is the effective Hamiltonian $s_T(p)=\overline H(p)$, see Fig. \ref{average-plane},  associated with the G-equation (a first or second order HJE), defined in a nonlinear eigenvalue problem (a.k.a cell or corrector problem \cite{LPV,E1989}):
\be\label{eq:generalGcell}
H(D^2v,p+Dv,x)=(\mathrm{or} \approx) \overline H(p) \quad \text{in $\Rset^n$}
\ee
where $v$ is called a corrector (eigenfunction), and $\bar H$ the eigenvalue.
Here $\approx$ refers to approximate solutions when exact solutions do not exist. See section 2 for more details. 

In general, the speed $\overline H(p)$ is direction $p$ dependent or anisotropic. Moreover,  except in some simple cases (e.g., 2D shear flows), there is no closed form formula of $\overline H$. If $\overline H(p)$ exists,  we have:
\be\label{eq:large-time-limit}
\overline H(p)=-\lim_{t\to \infty}{G(x,t) \over t},
\ee
where $G$ is the solution to the corresponding equation (\ref{main-eq}) subject to $G(x,0)=p\cdot x$.  The above formula can be used to numerically compute $\overline H(p)$, see 
\cite{Q2003,KLX_21,LXY_13a}.

\medskip

 {\bf $\bullet$ Perspective 2:} In combustion literature and experiments, the ratio $s_T/s_l$ is  measured by the ratio between the surface area of the wrinkled flame front and the flat one (unwrinkled).   According to (11) and (12) in \cite{KAW1988}, the turbulent flame speed at time $t$ is approximately the volume average:
$$
u_T(t)=\int_{[0,1]^n}s_l|DG(x,t)|\,dx.
$$
when the flow field $V$ is 1-periodic, incompressible with mean zero, 
see also \cite{CTVV2003} for relevant discussions.  
Integrating (\ref{main-eq}) over $x$ gives:
$$
\int_{[0,1]^n}s_l|DG(x,t)|\,dx=-\int_{[0,1]^n}G_t(x,t)\,dx.
$$
Then the time average of $u_T$ is 
$$
{1\over t}\int_{0}^{t}u_T(s)\,ds= -{\int_{[0,1]^n} \, (G(x,t) - p\cdot x ) \, dx \over t} {\sim \atop t\gg 1} -{\int_{[0,1]^n} \, G(x,t)  \, dx \over t}.
$$
As noted in \cite{KAW1988}, an open problem is whether the time average of $u_T$ 
converges to a constant as $t \uparrow \infty$, or the validity of the spatially averaged version of 
(\ref{eq:large-time-limit}). A computational example 
in \cite{KAW1988}
based on solutions of Navier-Stokes equations for the flow velocity $V$ showed the convergence of the time-averaged $u_T$. 

\begin{figure}
\begin{center}
\includegraphics[scale=0.8]{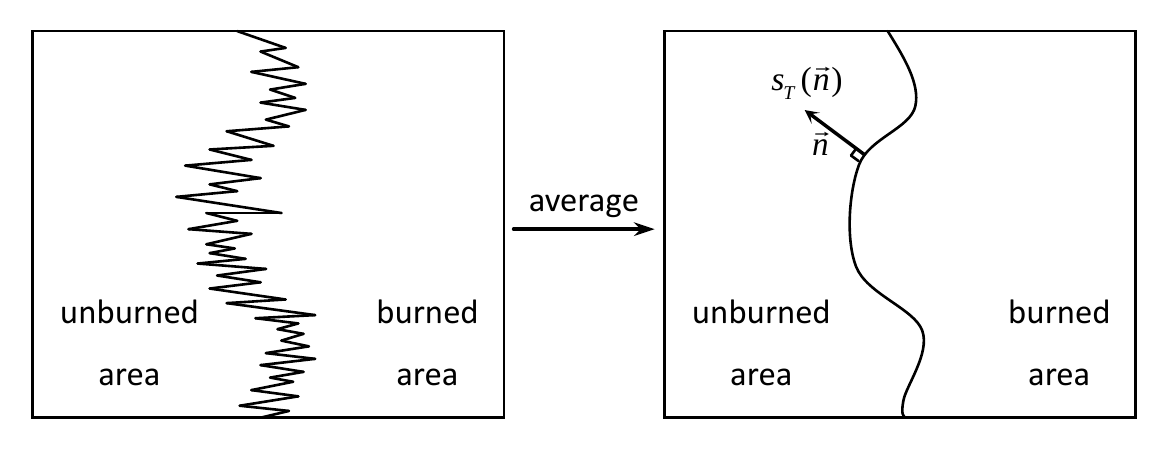}
\caption{Averaging fluctuations around a general flame front.}
\label{average-front-1}
\end{center}
\end{figure}

The existence of $\overline H(p)$ usually leads to the following full homogenization result in the single small scale situation \cite{E1992}. Denote by $\epsilon$ the Kolmogorov scale.  Change variables
$$
x\to {x\over \epsilon}, \quad V(x)\to V\left({x\over \epsilon}\right) \quad \mathrm{and} \quad d\to d\, \epsilon
$$
and denote by  $G^{\epsilon}(x,t)$ the corresponding solution to (\ref{main-eq}) with the prescribed initial flame front $G(x,0)=g(x)$. Then
$$
\lim_{\epsilon\to 0}G^{\epsilon}(x,t)=\bar G(x,t),
$$
where $\bar G(x,t)$ solves an effective equation
$$
\begin{cases}
\bar G_t+\overline H(D\bar G)=0\\[3mm]
\bar G(x,0)=g(x).
\end{cases}
$$
Following standard notations, we write $s_T(p)=\overline H(p)$. The effective equation says that the averaged flame front will propagate with speed
$$
v_{\vec{n}}=\overline H(\vec{n})
$$
along its normal direction $\vec{n}$. See Fig. \ref{average-front-1} for an illustration.

\medskip

We shall focus on the following two topics closely related to a central theme in combustion ({\it how fast can it burn ?} ).

\medskip

(1) Proving the existence of $\overline H(p)$ under physically meaningful assumptions. This is equivalent to finding solutions or approximate solutions of the form (\ref{steady-form}), which reduces to studying cell problems in homogenization theory of HJE \cite{LPV,E1989}. Due to the lack of coercivity from a large flow, standard PDE methods are insufficient.  The Lagrangian approaches are introduced, especially game theoretic methods in the presence of non-convexity from curvature and flame stretching effects. 

\medskip

(2) 
Studying qualitative and quantitative properties of  $\overline H(p)$ beyond existence, a mathematical topic of physical impact. Below are two issues addressed in a large body of combustion literature through theoretical or empirical approaches. 

\medskip

I. Dependence of the turbulent burning velocity on flow intensity. Precisely speaking,  scale the flow field $V$ to $AV$ by a positive constant $A$ and inquire about the growth behavior of  $\overline H(p)$ as $A\to +\infty$. This is one of the major approaches to speed up flame propagation within a combustion vessel. 

\medskip

II.  Understanding how the turbulent burning velocity relies on the curvature and strain effect.  In particular, an important problem in turbulent combustion is whether adding the curvature correction will decrease the prediction of the turbulent burning velocity. This is equivalent to asking whether $\overline H$ is a monotonically decreasing function of the Markstein length (or Markstein number) $d$ in (\ref{eq:generalspeed}) that mainly depends on the type of burning fuel. 

\medskip

   Tackling I and II often calls for challenging analysis of the underlying dynamical systems (e.g. on flow geometry and asymptotics), which can be made precise without making general assumptions by studying specific flows with physical origin.
   The two common candidates to start with are shear flows and two dimensional (2D) cellular flows.  Even for 2D shear flow where the PDE under study basically reduces to ODE, the associated mathematical questions could be highly non-trivial when curvature terms are involved. In the next step to analyze 3D flows, we choose a well-known example, the Arnold-Beltrami-Childress (ABC) flow \cite{Arnold_65,DFGHMS,FGM_93, CG95}, a steady solution to the 3D Euler equation with chaotic yet non-ergodic streamlines.  The integrable case of the ABC flow is the 2D cellular flow (a.k.a. BC flow).  
   Recent progress in constructing ballistic orbits \cite{XYZ,MXYZ_16,KLX_21}
   is necessary to address issue I for 
   the inviscid G-equation in the ABC flow.
\medskip

For curvature and strain G-equations, 
our integrated Lagrangian and Eulerian methodology established first for the inviscid G-equation continues to shine in a two-player game framework in the non-convex setting. Even though the classical Euler-Lagrange duality via Legendre transform is missing, the success of our blended Lagrangian-Eulerian analysis points to a {\it subtle and hidden duality} for future research to explore in geometric PDEs and flows. On the other hand, for the curvature G-equations in three and higher dimensions, the averaging (homogenization) breaks down in large enough intensity shear flows \cite{MMTXY2023}, indicating a non-convex phenomenon that an existence proof for the effective burning speed through volume preserving (incompressible) flows is dimension and flow dependent.  

   \medskip

   {\bf $\bullet$ Outline of the paper.} In section 2, we review the basics of viscosity solutions and homogenization theory. In sections 3-4, we  discuss different types of G-equations, present theoretical and numerical results for the basic and curvature G-equations based on control and two-player game representations as well as PDE methods. Section 5 gives a two-player differential game and PDE analysis for the strain G-equation.
    Section 6 revisits other passive scalar flame propagation models related to G-equation. Section 7 is a brief overview of stochastic homogenization and multi-scale modeling of 
   the basic G-equation in turbulent combustion. 
   Conclusions are in section 8. 
   The open problems appear along the way.

\section{Preliminary:  Viscosity Solutions and Homogenization} For the reader's convenience, we will review some basic concepts and techniques about viscosity solutions and homogenization theory most relevant to problems discussed in this paper.  See \cite{CIL_92, Evanbk, T2021, ES2008}  for more details.
\setcounter{equation}{0}

Consider a general second order Hamilton-Jacobi equation (HJE)
\be\label{eq:2ndorder}
u_t+H(D^2u, Du, u, x)=0 \quad \text{on $\Rset^n\times (0,\infty).$}
\ee
Here the Hamiltonian $H=H(M,p,z, x)$ is a function defined on $(M,p,z, x)\in S^{n\times n}\times \Rset^n\times \Rset\times  \Rset^n$, where $S^{n\times n}$ is the set of $n\times n$ symmetric matrices.  We say that the Hamiltonian $H$
is {\bf coercive} in the gradient variable $p$ if for any fixed $K\geq  0$, 
$$
\lim_{|p|\to +\infty}\min_{|M|+|z|+|x|\leq K}H(M,p,z, x)=+\infty.
$$

\begin{dfn} An upper-semicontinuous function $u(x,t)$ defined on $\Rset^n\times (0,\infty)$ is called a viscosity subsolution of (\ref{eq:2ndorder}) if for any $\phi\in C^2(\Rset^n\times (0,\infty))$ and any point $(x_0,t_0)\in \Rset^n\times (0,\infty)$, when
$$
0=\phi(x_0,t_0)-u(x_0,t_0)\leq \phi(x,t)-u(x,t)
$$
holds for $(x,t)$ in a neighborhood  of $(x_0,t_0)$,  
we have 
$$
\phi_t(x_0,t_0)+H(D^2\phi(x_0,t_0), D\phi(x_0,t_0), \phi(x_0,t_0), x_0)\leq 0. 
$$
\end{dfn}

\begin{dfn} An lower-semicontinuous function $u(x,t)$ defined on $\Rset^n\times (0,\infty)$ is called a viscosity supersolution of (\ref{eq:2ndorder}) if for any $\phi\in C^2(\Rset^n\times (0,\infty))$ and any point $(x_0,t_0)\in \Rset^n\times (0,\infty)$, when 
$$
0=\phi(x_0,t_0)-u(x_0,t_0)\geq \phi(x,t)-u(x,t) 
$$ 
holds for $(x,t)$ in a neighborhood  of $(x_0,t_0)$, we have 
$$
\phi_t(x_0,t_0)+H(D^2\phi(x_0,t_0), D\phi(x_0,t_0), \phi(x_0,t_0), x_0)\geq 0. 
$$
\end{dfn}

\begin{dfn} A continuous function $u(x,t)$ defined on $\Rset^n\times (0,\infty)$ is called a viscosity solution of (\ref{eq:2ndorder}) if it is both a viscosity subsolution and a viscosity supersolution. 
    
\end{dfn}

Similarly, we can define viscosity subsolutions, supersolutions and solutions  for steady state equations 
\be\label{eq:2dsteady}
H(D^2u, Du, u, x)=0 \quad \text{on $U$},
\ee
where $U$ is an open subset of $\Rset^n$. 

Moreover, we say that (\ref{eq:2dsteady}) has {\it approximation solutions} if for any $\delta>0$, there is a continuous function $u_\delta(x,t)$ that satisfies the following inequality in the viscosity sense
\be\label{eq:approxs}
-\delta\leq H(D^2u^{\delta}, Du^\delta, u^\delta, x)\leq \delta  \quad \text{on $\Rset^n\times (0,\infty)$}.
\ee

The notion of viscosity solutions provides a rigorous mathematical framework to describe the correct ``physical" solution of  the corresponding equations (\ref{eq:2ndorder} ) or (\ref{eq:2dsteady}) when classical solutions might not exist. Important examples include first order HJ equations or second order degenerate elliptic equations arising from control theory or front propagation problems in applications.  Heuristically, the name comes from the vanishing viscosity method.  To stabilize the equation numerically, 
an artificial viscosity term is added to the equation although often without clear physical meanings. Precisely speaking, given $ h >0$, let $u^{h}$ be a smooth solution (if it exists) to 
$$
-h\, \Delta u^{h}+u_t^{h}+H(D^2u^h, Du^h, u^h, x)=0 \quad \text{on $\Rset^n\times (0,\infty)$}
$$
with  given initial data $u^h(x,0)=g(x)$.  Then it is natural to expect that $u^{h}$ converges to the correct physical solution of (\ref{eq:2ndorder}) as $h\to 0$ (vanishing viscosity). In fact, using maximum principle, it is not hard to prove that if 
$$
\lim_{h\to 0}u^{h}(x,t)=u(x,t) \quad \text{locally uniformly in $\Rset^n\times (0,\infty)$},
$$
then $u$ is the viscosity solution to (\ref{eq:2ndorder}) subject to $u(x,0)=g(x)$.  Similar limiting process can also be carried out for the steady case (\ref{eq:2dsteady}) with given boundary data $u$ on $\partial U$.  Although this vanishing viscosity limit is intuitively clear, the existence of viscosity solutions is usually established by Perron's method that is more convenient.  Under suitable assumptions, uniqueness of viscosity solutions holds  with prescribed initial data or boundary data. We refer to \cite{CIL_92} for more details. 

Next we review the basics of homogenization of nonlinear equations that started from \cite{LPV}. In our context,  assume that $H=H(M,p,x)$ is periodic in $x$ and consider the solution $u^{\epsilon}$ to the oscillatory equation
\be\label{eq:2ndorderh}
\begin{cases}
u_{t}^{\epsilon}+H\left(\epsilon D^2u^{\epsilon}, Du^{\epsilon}, {x\over \epsilon}\right)=0\qquad \text{in $\Rset^n\times (0,\infty)$}\\
u^{\epsilon}(x,0)=g(x).
\end{cases}
\ee
when the underlying environment has small scale $\epsilon>0$.  Here we omit the dependence on the $u$ variable for convenience. Under proper assumptions, it can be shown \cite{LPV, E1992} that  
\be\label{eq:h-convergence}
\lim_{\epsilon\to 0}u^{\epsilon}(x,t)=u(x,t) \quad \text{locally uniformly in $\Rset^n\times (0,\infty)$}.
\ee
 Here $u(x,t)$ is the unique continuous viscosity solution to the effective equation
\be\label{eq:effective1}
\begin{cases}
\bar u_t+\overline H(D\bar u)=0\\
\bar u(x,0)=g(x).
\end{cases}
\ee
The above limiting process is called  homogenization.
Given $p\in \Rset^n$,  the effective Hamiltonian $\overline H(p): \Rset^n\to \Rset$ is the unique number such that the following cell (corrector) problem
\be\label{eq:cellp}
H(D^2v, p+Dv, x)=\overline H(p) \quad \text{in $\Rset^n$}
\ee
has a continuous $\Zset^n$-periodic viscosity solution (or approximate solutions in the sense of (\ref{eq:approxs})) $v(x)$, which is called a corrector.  Here we would like to stress that the solution $v$ in general might not be unique even up to an additive constant.  When the original Hamiltonian $H$ is from front propagation problems, the effective Hamiltonian $\overline H$ has a clear physical meaning of effective propagation speeds in the corresponding setting. For the mechanical Hamiltonian $H(p,x)={1\over 2}|p|+V(x)$, the effective Hamiltonian and the associated cell problem have helped construct quasi-modes of Sch$\Ddot{\mathrm{o}}$rdinger operators via various quantization procedures \cite{E2004}. 

If the existence of the cell problem can be established (i.e., existence of effective Hamiltonian), then the convergence (\ref{eq:h-convergence}) usually follows from the perturbed test function method \cite{E1989,E1992}. Moreover,  comparison principle implies that if the initial value $u(x,0)=p\cdot x$, then
\be\label{eq:l-time-limit}
-\lim_{t \to +\infty}{u(x,t)\over t}=\overline H(p) \quad \text{for all $x\in \Rset^n$}.
\ee

 Technically speaking, the key part is to establish the existence of the cell problem (\ref{eq:cellp}). Below is the mechanism first introduced in \cite{LPV,E1992}: for $\lambda>0$, let $v_\lambda$ be the unique periodic viscosity solution of 
$$
\lambda \,  v_{\lambda}+H(D^2v_\lambda, p+Dv_\lambda, x)=0  \quad \text{in $\Rset^n$}.
$$
Then show that 
$$
-\lim_{\lambda\to 0}\lambda \, v_{\lambda}(x)=\text{ a constant} \quad \text{for $x\in \Rset^n$}.
$$
This constant, if it exists, is the effective Hamiltonian $\overline H(p)$.  To prove the above, it suffices to show that the oscillation of $v_\lambda$ is uniformly bounded:
\be\label{eq:uniformbound}
\max_{x,y\in \Rset^n}|v_{\lambda}(x)-v_{\lambda}(y)|\leq C
\ee
for a constant $C$ independent of $\lambda$.  The standard approach to prove the uniform bound is to derive equi-continuity of $v_\lambda(x)$ for $\lambda\in [0,1]$ (e.g. uniform H\"older continuity 
\cite{LPV}). Then $v=\lim_{\lambda\to 0}(v_\lambda(x)-v_\lambda(0))$ (up to a subsequence if necessary) is a solution to the cell problem. 

  However, when $H$ is not coercive in the gradient variable $p$ or is only degenerate elliptic in the Hessian variable $M$, like the curvature G-equation,  the equi-continuity is usually not valid. There appears no systematic way to obtain (\ref{eq:uniformbound}). Customized methods based on the structure of equations are often needed, often relying on approximate solutions of the cell problem (\ref{eq:cellp}) in the sense of (\ref{eq:approxs}). 

  From numerical point of view \cite{ES2008}, the effective equation can be used to approximate the original oscillatory equation (\ref{eq:2ndorderh}) where the computation is very expensive in order to resolve the small scale $\epsilon$. Of course, one has to first compute the effective Hamiltonian $\overline H(p)$, which can be viewed as some sort of nonlinear averaging of the original Hamiltonian. However, even for simple cases like $H(p,x)={1\over 2}|p|^2+V(x)$ or $H(p,x)=a(x)|p|$,  the effective Hamiltonian does not have close-form formulas except in 1D. Various numerical schemes \cite{Q2003,LXY_13a,KLX_21} have been developed to compute $\overline H(p)$ based on the large time limit (\ref{eq:l-time-limit}). 

   Understanding nontrivial analytic  properties of the effective Hamiltonian is extremely challenging and  often requires deep tools beyond standard PDE techniques. Such kinds of problems have been studied for decades in the areas of first passage percolation in probability theory, stable norm/$\beta$-functions in geometry, and dynamical systems, which are essentially effective Hamiltonian in the corresponding context. We refer to \cite{TY2022} and reference therein for more discussion about such connections. Nevertheless, there are not much relevant works in the PDE literature. See \cite{JTY2021,TY2023} for recent progress in this direction. In particular, it is surprising  to see that local properties of $\overline H$ could be dramatically different from the original Hamiltonian. For example, when $n=2$, for generic potential function  $V$,  the effective Hamiltonian $\overline H(p)$ associated with $H(p,x)={1\over 2}|p|^2+V(x)$ is piece-wise 1D for $p$ in a dense open set $U\subset \Rset^2$ (\cite{Y2022}).

   Finally,  we would like to mention that there is a large body of literature and active research on homogenization of scalar conservation laws and Hamilton-Jacobi PDEs in random environments, see \cite{WX_97,S1999, RT2000} for early works. It is an interesting question to explore how stochastic modeling of turbulent flows $V$ in practical situations is related to typical mathematical assumptions in stochastic homogenization. In this survey paper, we will mainly focus on the periodic setting, with 
   the stochastic aspects of G-equation briefly reviewed in section 7.

   \section{Inviscid G-equation: the Basic Case}
   \setcounter{equation}{0}
Let $G(x,t)$ be the unique viscosity solution to the inviscid G-equation
\be\label{in-eq}
\begin{cases}
G_t+|DG|+V(x)\cdot DG=0 \quad \text{on $\Rset^n\times (0,\infty)$}, \\
G(x,0)=g(x).
\end{cases}
\ee
As a solution to a convex HJE, $G(x,t)$ has a control formula \cite{Evanbk}:
\be
G(x,t)=\inf_{\alpha\in \mathcal{A}_t}g(\xi_\alpha(t)). \label{controlform}
\ee
Here $\mathcal{A}_t$ is the set of measurable functions $\alpha=\alpha(s): [0,t]\to \overline {B_1(0)}$  and $\xi_\alpha$ is the Lipschitz continuous curve satisfying
$$
\begin{cases}
\dot \xi_\alpha(s)=-V(\xi_\alpha(s))+\alpha(s)\quad \text{for a.e. $s\in [0,t]$}\\
\xi_\alpha(0)=x.
\end{cases}
$$
See \cite{Evanbk} for connections between general convex Hamilton-Jacobi equations and control theory. We define the following reachability notion. 

\begin{dfn} Given $x,y\in \Rset^n$, we say that $y$ is reachable  from $x$ within time $T$ if there exists  a Lipchitz continuous curve $\xi:[0,T]\to \Rset^n$ such that $\xi(0)=x$, $\xi(T)=y$ and $|\dot \xi(s)+V(\xi(s))|\leq 1$ for a.e. $s\in [0,T]$. 
\end{dfn}

Then we have the following one-sided bound. 

\begin{lem}\label{one-side-bound1} Suppose that $u$ is a viscosity subsolution to 
$$
|Du|+V(x)\cdot Du\leq C  \quad \text{in $\Rset^n$}
$$
for some fixed constant $C$. If $y$ is reachable from $x$ within time $T$, then
$$
u(x)\leq u(y) + CT. 
$$
\end{lem}
Proof: For convenience,  we may assume that $u$  and the curve $\xi$ connecting $x$ and $y$  are both $C^1$. Then
$$
{du(\xi(t))\over dt}=Du\cdot \dot \xi\geq -V(\xi)\cdot Du-|Du|\geq -C. 
$$
Hence
$$
u(y)-u(x)=\int_{0}^{T}{du(\xi(t))\over dt}\geq -CT. 
$$
\qed
\subsection{Homogenization of invicid G-equation}
 For convenience,  we assume that the flow field  $V$ is $\Zset^n$-periodic, $C^1$ and incompressible (i.e., $\mathrm{div}(V)=0$). The incompressibility assumption can be relaxed.  Meanwhile, it is easy to give examples of smooth $V$ where homogenization fails. The following  homogenization  result of inviscid G-equation was proved independently in \cite{XY_10} and \cite{CNS} through distinct methods.  For $\epsilon>0$, let $G^{\epsilon}(x,t)$ be the solution to 
 \be\label{invG}
\begin{cases} 
 G_{t}^{\epsilon}+|DG^{\epsilon}|+V\left({x\over \epsilon}\right)\cdot DG^{\epsilon}=0 \quad \text{on $\Rset^n\times (0,\infty)$} \\
 G^{\epsilon}(x,0)=g(x).
\end{cases}
\ee

\begin{thm} \label{mainthm}
$$
\lim_{\epsilon\to 0}G^{\epsilon}(x,t)=\bar G(x,t) \quad \text{locally uniformly on $\Rset^n\times (0,\infty)$}.
$$
Here $\bar G$ is the unique viscosity solution to the effective equation
\be
\begin{cases}
 {\bar G}_t+\overline H(D\bar G)=0 \quad \text{on $\Rset^n\times (0,\infty)$},\\
\bar G(x,0)=g(x). \label{aveGeq}
\end{cases}
\ee
$\overline H(p)\in C(\Rset^n, [0,\infty))$ is a convex positive homogeneous function of degree one and satisfies the enhancement property, i.e.,  $\overline H(p)\geq |p|$.

\end{thm}

\subsubsection{Sketch of proof.} As mentioned in the previous section, the key is to prove the existence of periodic solution or approximate solution to the cell problem for any given vector $p\in \Rset^n$:
$$
|p+Dv|+V(x)\cdot (p+Dv)=\overline H(p) \quad \text{in $\Rset^n$}.
$$
Below we present the method in \cite{XY_10} with a slightly refined version in \cite{XY_14a}.
\medskip

{\bf Step 1}:  Consider the $\lambda$-auxiliary equation for $\lambda>0$:
$$
\lambda \, v_{\lambda}+|\, p+Dv_{\lambda}|+V(x)\cdot (\;  p+Dv_{\lambda})=0.
$$
The goal is to show that 
$$
\lim_{\lambda\to 0}\lambda\, v_{\lambda}(x)=\mathrm{constant}  \quad \text{uniformly on $\Rset^n$}.
$$
This will follow immediately if we can establish the uniform bounded-ness of the oscillation: 
$$
\max_{x,y\in \Rset^n}|v_{\lambda}(x)-v_{\lambda}(y)|\leq C
$$
for a constant $C$ independent of $\lambda$.  Due to the lack of coercivity, the standard uniform Lipschitz continuity estimate of $v_\lambda$ is not available.  To overcome that, we consider
$$
\bar v(x)=\limsup_{y\to x,\ \lambda\to 0}\lambda\;  v_{\lambda}(y).
$$
Then routine argument shows that $\bar v$ is an upper-semicontinuous solution to
$$
|D\bar v|+V(x)\cdot D\bar v\leq 0. 
$$
Integrating over the unit cube $Q_1=[0,1]^n$ on both sides leads to 
$$
\int_{Q_1}|D\bar v|\,dx=0.
$$
Accordingly, $|D\bar v|\equiv 0$ and $\bar v(x)=\mathrm{constant}$ that is denoted by $v^*$.

\medskip

{\bf Step 2:}  For $x\in \Rset^n$, set
$$
\underline{v}(x)=\liminf_{y\to x, \ \lambda\to 0}\lambda \; v_{\lambda}(y).
$$
Due to the convexity of the Hamiltonian $H(p,x)=|p|+V(x)\cdot p$ with respect to $p$, $h_\lambda(x)=-v_\lambda$ is a viscosity subsolution to 
$$
-\lambda \, h_{\lambda}+|\, p-Dh_{\lambda}|+V(x)\cdot (\;  p-Dh_{\lambda})\leq 0.
$$
Note that this is in general not true for viscosity solutions of  non-convex HJ equations.  Similar to Step 1, $\limsup_{\lambda \to 0}\, \lambda \, h_\lambda(x)=-\underline{v}(x)$ is a constant. Hence $\underline v(x)=\mathrm{constant}$ that is denoted by $v_*$.

Obviously, $v_*\leq v^{*}$.  It suffices to prove that
$$
v_*\geq v^{*}
$$
By choosing a suitable subsequence of $\lambda$, this follows from the next step and Lemma \ref{one-side-bound1}.

\medskip

{\bf Step 3:} It is easy to see that for any $x\in \Rset^n$, there exist $y\in \Rset^n$, $r>0$ and $T>0$ such that for every pair $(x',y')\in B_r(x)\times B_r(y)$,  $y'$ is reachable from $x'$ within time $T$.

\begin{rmk} The above method and conclusion can be easily extended to the time dependent case where $V=V(x,t)$ and is periodic in $(x,t)$. The key is to recover step 1, i.e. showing that the limsup is a constant. The other parts are similar.  In fact, assume that $\bar v(x,t)$ is periodic in both time and space and a viscosity subsolution to 
$$
\bar v_t+|D\bar v|+V(x,t)\cdot D\bar v\leq 0 \quad \text{on $\Rset^n\times \Rset$}.
$$
As in the previous step 1, integrating with respect to both $x$ and $t$ on $Q_1\times [0,1]$ leads to $D\bar v=0$, i.e., $\bar v$ is constant in $x$.   Then $v_t\leq 0$, which implies that $\bar v$ is also constant in $t$ due to the periodicity in $t$.  
\end{rmk}

\subsubsection{Optimal Convergence Rate $O(\epsilon)$} 
We expect to obtain the optimal convergence rate $O(\epsilon)$  if the initial data $g$ is Lipschitz continuous:
\be\label{optimal-rate}
|G^{\epsilon}(x,t)-\bar G(x,t)|\leq O(\epsilon). 
\ee
This $O(\epsilon)$ optimal convergence rate has been established in \cite{TY2022} for general coercive convex Hamilton-Jacobi equations. See \cite{TY2022} and references therein for other important works in this topic.  

Below we list the main steps with details left to the interested readers to explore. It is enough to prove this for $(x,t)=(0,1)$. 

\medskip

Step 1: For $x,y\in \Rset^n$, define a distance function
$$
d(x,y)=\text{the minimum time to reach $y$ from $x$}. 
$$
See \cite{NN_11, CS_13} for some interesting estimates about $d(x,y)$ although we only need $d(x,y)<\infty$  for our purpose. Note that $d$ might not be continuous due to the lack of coercivity.
\medskip

Step 2: This is the key step.  Owing to \cite{B1992} and the corresponding connections elaborated in \cite{TY2022},  there exists a continuous function $\bar d(x,y)$ such that
$$
\left|\epsilon d\left(0, {x\over \epsilon} \right)- \bar d(0,x)\right|\leq C\epsilon |x| 
$$
for a constant $C$ independent of $\epsilon$. We want to point out that to establish the upper bound $\epsilon d\left(0, {x\over \epsilon} \right)\leq \bar d(0,x)+C\epsilon|x|$ is not very hard. The main difficulty is to verify the other direction that relies on a surprising and  beautiful cutting lemma (Lemma 2 in \cite{B1992}).  Then (\ref{optimal-rate}) should follow by establishing connection between (\ref{controlform}) and the control formulation of the solution to the effective equation (\ref{aveGeq}) in \cite{NN_11}.


\subsection{Explicit formula of $\overline H(p)$ in 2D shear flows}
For 2D shear flow $V$ (e.g., the flow within a slot burner), we obtain the explicit formula of $\overline H(p)$.
Let $V(x_1,x_2)=(0,f(x_1))$ for a continuous periodic function $f:\Rset\to \Rset$. Given $p=(a,b)$, the corresponding cell problem reduces to the ODE:
$$
\sqrt{(a+v'(y))^2+b^2}+b\, f(y)=\overline H(p),
$$
where $v:\Rset\to \Rset$ is a periodic Lipschitz continuous function. 
Then
$$
|a+v'(y)|=\sqrt{(\overline H(p)-b\, f(y))^2-b^2}.
$$
It is not hard to derive that
\medskip

(1) when $|a|\leq \int_{0}^{1}\sqrt{(M-bf(y))^2-b^2}\, dy$,
$$
\overline H(p)=M=|b|+\max_{\Rset}\, b\, f;
$$

(2) when $|a|>\int_{0}^{1}\sqrt{(M-b\, f(y))^2-b^2}\, dy$, $\overline H$ is given by the following implicit formula
$$
|a|=\int_{0}^{1}\sqrt{(\overline H(p)-b\, f(y))^2-b^2}\, dy.
$$
The computation is similar to that in \cite{LPV} for $H(p,x)=|p|^2+W(x)$. The detail is left to interested readers.

\subsection{Growth law of $\overline H(p)$ in 2D cellular flows}

In general, the behavior of $\overline H(p)$ as a function of $p$ is quite complicated and does not possess any  formula in closed form.  Here we focus on how $\overline H(p)$ depends on the large flow intensity, an extensively studied issue in combustion literature. Precisely speaking, for $A>0$, let $V(x)\to AV(x)$ and consider the corresponding $\overline H(p,A)$ in the cell problem:
$$
|p+Dv|+AV(x)\cdot (p+Dv)=\overline H(p,A) \quad \text{in $\Rset^n$}.
$$

\begin{ques}
What is the growth pattern of $\overline H(p,A)$ as $A\to +\infty$ ?
\end{ques}

\medskip

To derive meaningful result,  we need to look at physically interesting examples of $V$. The first is a 2D Hamiltonian flow, the cellular flow with a typical form:
$$
V(x_1,x_2)=(-H_{x_2},\  H_{x_1}).
$$
Here the stream function $H=\sin x_1\sin x_2$.  As mentioned before, this is an integrable form (after suitable rotation) of the ABC flow \cite{Arnold_65,DFGHMS,FGM_93},  a steady periodic solution to the 3D Euler equation. Using (\ref{eq:l-time-limit}) and the control formulation (\ref{controlform}),  it was proved in \cite{XY_13} that for any unit vector $p\in \Rset^2$
\be\label{Ggrowth}
{A\, \pi(|p_1|+|p_2|)\over 2\log A+C_2}\leq \overline H(p,A) \leq {A\, \pi(|p_1|+|p_2|)\over 2\log A+C_1}.
\ee
Here $C_1$ and $C_2$ are two constants independent of $A$ and $p$.  Below is an interesting question:
\medskip

 \begin{ques}
Does there exist a constant $C$ such that
$$
\left|\overline H(p,A)-{A\, \pi(|p_1|+|p_2|)\over 2\log A+C}\right|\leq O(A^{\alpha})
$$
for some $\alpha\in [0,1)$?
\end{ques}
This question is closely related to identifying the curvature effect discussed in a later section.  

  In numerical computation  \cite{KAW1988} and physical modeling \cite{Clavin_85,FIL_97} of 
  G-equation, a viscosity term is often added (explicitly or implicitly), which leads to the viscous G-equation:
  $$
  -d\Delta G+|DG|+AV(x)\cdot DG=0.
  $$
Accordingly, a natural question is to identify the effect of viscosity term in lieu of the curvature term.  Let $\overline H(d,p,A)$ be the unique constant such that the  corresponding cell problem
$$
-d\Delta v+ |p+Dv|+AV(x)\cdot (p+Dv)=\overline H(d, p,A) \quad \text{in $\Rset^n$}
$$
has a periodic solution $v$. Due to the presence of the $\Delta v$,  the existence of $\overline H(d,p,A)$ is well-known (see \cite{E1992} for instance). Surprisingly, it turns out that the viscosity term will dramatically slow down the effective burning velocity. Precisely speaking, it was proved in \cite{LXY_11} that
\begin{thm} For fixed $d>0$ and $V(x_1,x_2)=(-H_{x_2},\  H_{x_1})$ with $H=\sin x_1\sin x_2$,
$$
\sup_{A\geq 0}\overline H(d, p,A)\leq C_d\, |p|
$$
for a constant $C_d$ depending only on $d$. 
\end{thm}
 Below is the outline of the proof.  Integrating  both sides of the above cell problem gives 
$$
\overline H(d, p,A)=\int_{\Bbb T^n}|p+Dv|\,dx.
$$
Accordingly, we need to derive a uniform $L^1$ bound of $Dv$.  This is done in two steps:
\medskip

Step 1: Consider periodic solution  $T$ of the linear equation
$$
-d\Delta T+AV(x)\cdot (p+DT)=0 \quad \text{in $\Rset^n$}.
$$
It was first established that $L^1$ norms of $D\, v$ and $D\, T$ are comparable:
$$
{1\over C}||D\, T||_{L^1(\Bbb T^2)}\leq ||D\, v||_{L^1(\Bbb T^2)}\leq C||D\, T||_{L^1(\Bbb T^2)}
$$
for a constant $C$ independent of $A$. 

\medskip

Step 2:  Prove that the $L^1$ norm of $Dv$ is uniformly bounded
$$
||DT||_{L^1(\Bbb T^2)}\leq C.
$$
This is achieved by (1) carefully using interior decay estimates of $||DT||_{L^2}$ established in \cite{NPR2005} and (2) delicate estimates of $||DT||_{L^1}$ near the boundary $\{H=0\}$, especially near corners, via suitable change of variables. 

\subsection{Growth law of $\overline H(p)$ in the ABC flow and Kolmogorov flow} In three dimensional (3D) incompressible flows, due to the presence of chaotic structures, it becomes much more challenging to analyze the dependence of $\overline H(p,A)$ on $A$.  The following asymptotic result was proved in \cite{XY_14b}:
\begin{thm}\label{thm:a-limit}
$$
\begin{array}{ll}
\lim_{A\to +\infty}s_T(p,A)/ A&=\max_{\sigma\in \Lambda}\int_{\Bbb T^n}p\cdot V(x)\,d\sigma\\[5mm]
&=\max_{\dot \xi=V(\xi)}\limsup_{T\to +\infty}\xi (T)\cdot p/ T.
\end{array}
$$
where $\Lambda$ is the collection of all Borel probability measures on $\Bbb T^n$ which are invariant under the flow $\dot \xi=V(\xi)$. 
\end{thm}

Below is a basic question.

\medskip

\begin{ques}
For a given 3D periodic incompressible flow, determine whether
\be\label{positive}
\lim_{A\to +\infty}{s_T(p,A)\over A}>0 \; 
\ee
or equivalently whether $\overline H(p,A)$ grows linearly with respect to $A$.
\end{ques}
This seemingly simple problem is  in general very challenging when $n\geq 3$ since it requires to identify the long time behavior of trajectories of $\dot \xi=V(\xi)$ that usually has chaotic structures. In particular, to prove the linear growth might need to find  unbounded trajectories with predictable long time behavior (e.g. 
 periodic or quasi-periodic when projected on the flat torus $\Bbb T^n$).  For example, if there is a periodic trajectory with a rotation vector not perpendicular to $p$, then (\ref{positive}) holds. 
 \medskip
 
Now let us turn to a well-known example of 3D incompressible flow introduced by Arnold, Beltrami and Childress, the so called ABC flow \cite{Arnold_65,DFGHMS,FGM_93}. For $x=(x_1,x_2,x_3)\in \Rset^3$,  
$$
V(x)=(A\sin x_3+C\cos x_2, B \sin x_1+A\cos x_3, \ C\sin x_2+B\cos x_1), 
$$
where $A$, $B$ and $C$ are three positive constants. Note that the  parameter $A$ here is not the flow intensity.  If one of them is zero, the flow is up to a $\pi/4$ rotation the integrable cellular flow.  The  ABC flow is a periodic steady solution of the Euler equation. In fact, the ABC flow satisfies the Beltrami property ($\nabla=({\partial \over \partial x_1},{\partial \over \partial x_2}, {\partial \over \partial x_3})$ denotes the gradient operator):
$$
\nabla\times V\, :=\mathrm{curl}(V)=V. 
$$
Using the following identity in vector calculus
$$
(V\cdot \nabla) V= \nabla\, {1\over 2}|V|^2-V\times (\nabla \times V),
$$
we have that the second term on the right hand side vanishes and 
$$
V_t+(V\cdot \nabla) V=(V\cdot \nabla) V=-\nabla P \quad  \mathrm{and} \quad  \nabla\cdot V=0,
$$
where the pressure $P=-{1\over 2}|V|^2$. Hence ABC flow is a steady state of the 3D incompressible Euler equation in fluid mechanics, see \cite{DFGHMS,FGM_93} for more discussion of its dynamic properties. 
\medskip

\begin{figure}
\begin{center}
\includegraphics[scale=0.3]{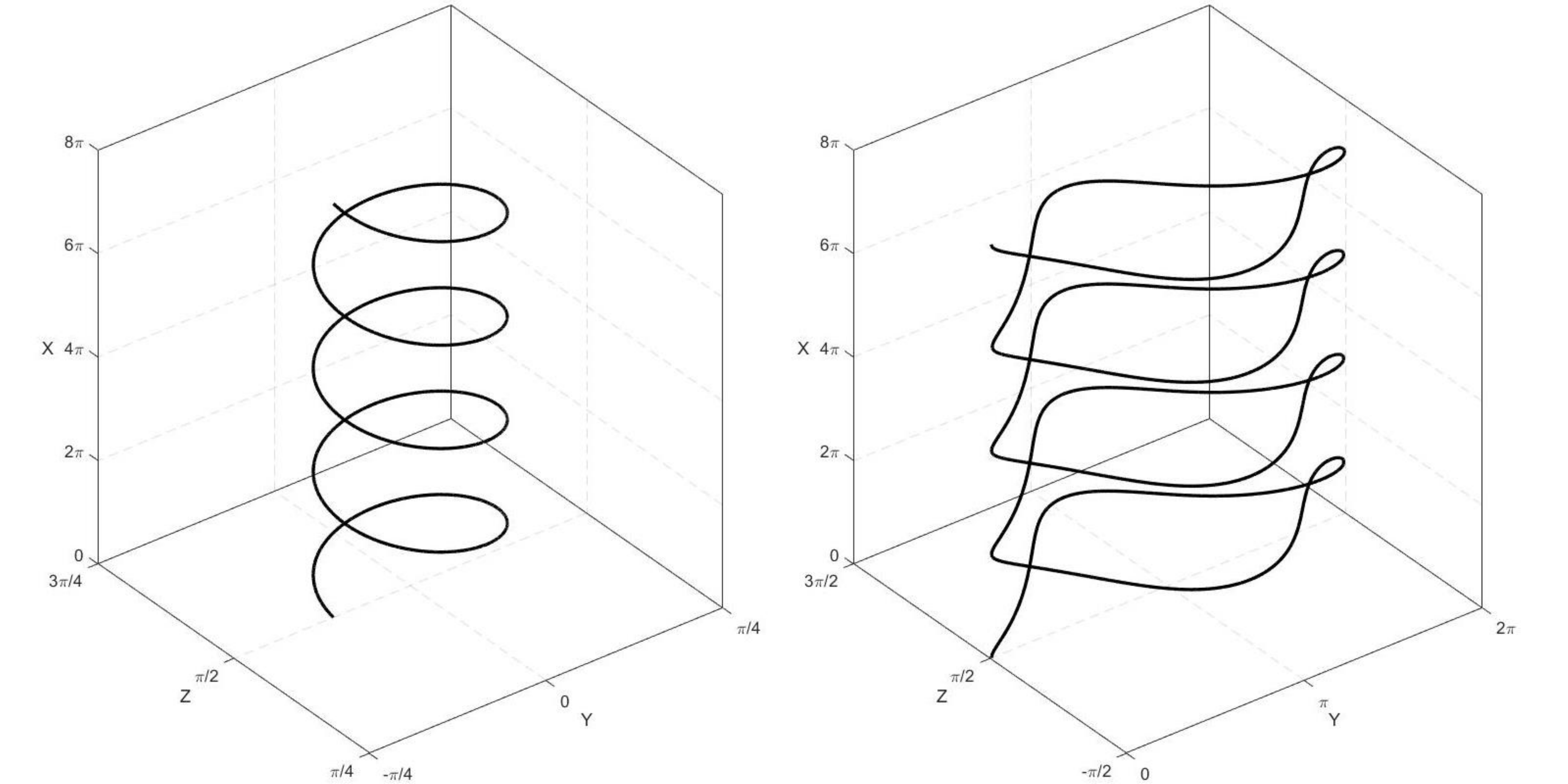}
\caption{Ballistic orbits, periodic (modulo $2\pi$) along $x$-direction, in 1-1-1 ABC (left) and Kolmogorov (right)
flows.}
\label{borb}
\end{center}
\end{figure}

 In the symmetric case $A=B=C=1$, the existence of non-contractable periodic orbits (in the sense stated precisely below) has been proved by using symmetry of the flow \cite{XYZ}:
 \begin{thm}
 There exists  $t_0>0$ and  a solution  $X(t)=(x(t),y(t), z(t))$ to the 1-1-1 ABC flow system such that for each $t\in\Rset$ we have
$$
X(t+t_0)=X(t)+(2\pi, 0, 0).
$$  
Then  $Y(t)=(z(t), x(t), y(t))$ and $Z(t)=(y(t),   z(t),  x(t))$ are clearly also solutions and satisfy 
$$
Y(t+t_0)=Y(t)+(0, 2\pi, 0)   \qquad \mathrm{and}   \qquad Z(t+t_0)=Z(t)+(0, 0, 2\pi).
$$
\end{thm}
The Kolmogorov flow with more chaotic phase space \cite{CG95} is:
\be
V(x, y, z) = (\sin z, \sin x, \sin y), 
\label{flow}
\ee
for which similar ballistic orbits exist \cite{KLX_21}, e.g. those that extend along each coordinate direction. See Fig. \ref{borb} for an illustration of two such orbits extending in a spiral manner in the $x$ direction. For related studies, we refer to \cite{MXYZ_16,KLX_21}. 
\medskip

Thanks to the control formula (\ref{controlform}), an immediate corollary is that in the 1-1-1 ABC and Kolmogorov flows,  the effective burning velocity grows linearly with respect to the flow intensity $A$, i.e., 
$$
\lim_{A\to +\infty}{\overline H(p,A)\over A}>0 \quad \text{for all unit vector $p\in \Rset^3$.}
$$
Moreover for $p$ along a coordinate direction \cite{KLX_21},
\be
{2\pi \over t_0}\, A + {2\pi \over t_0 \|V\|_{\infty}} 
\leq \overline H(p,A) \leq 
\|V \cdot p\|_{\infty}\, A + 1. \label{lin-growth}
\ee
With numerical estimates of period values $t_0$ from the 1-1-1 ABC and Kolmogorov flows, the following inequalities \cite{KLX_21} follow from (\ref{lin-growth}):
\[
1.942 A + 0.793 \leq \overline H \leq 2A + 1, \;\; (\text{1-1-1 ABC});
\]
\[
0.414 A + 0.239 \leq \overline H \leq A + 1, \;\; (\text{Kolmogorov}).
\]

\medskip

{\bf $\bullet$ Connection with the  Weinstein Conjecture} 

\medskip

To find non-contractible periodic orbits for general Beltrami flows is much more difficult and is closely related to the well known Weinstein conjecture due to the relation between Beltrami and Reeb vector fields (see \cite{MO2021} for instance).  Given a  Beltrami flow $V(x)$, we can define a one-form on the flat torus $\Bbb T^3$:
$$
a=V(x)\cdot dx. 
$$
Since $\mathrm{curl}(V)=\lambda V$ for  a constant $\lambda \not = 0$ (Beltrami property), 
$$
a\wedge da= \lambda \, |V(x)|^2dx_1\wedge dx_2\wedge dx_3.
$$
Accordingly,   ${V\over |V|^2}$ is a Reeb vector field if $V$ is nowhere vanishing.  Owing to a celebrated result of Taubes \cite{T2007},  $V$ has at least one unbounded periodic orbit with  certain rotation vector $q$, which implies that
$$
\lim_{A\to +\infty}{\overline H(p,A)\over A}>0, \quad \text{if  $p\cdot q\not=0$}.
$$
However, this is not enough to obtain linear growth along every unit direction. Note that Taubes' result does not apply to the 1-1-1 ABC flow that has critical points. 

\section{Curvature G-equation}
\setcounter{equation}{0}

Recall the curvature G-equation:
\be
G_t + \left(1- d\, \, \mathrm{div}\left({DG\over |DG|}\right)\right)_{+}|DG|+V(x)\cdot DG=0. \label{cG1}
\ee
A fundamental question is to understand how the curvature term influences the prediction of the effective burning velocity.  The first step is to rigorously establish the existence of effective burning velocity. Precisely speaking,  given a unit direction $p\in \Rset^n$, let $G(x,t)$ be the unique solution to 
\be\label{eq:linearcur}
\begin{cases}
G_t + \left(1- d\, \, \mathrm{div}\left({DG\over |DG|}\right)\right)_{+}|DG|+V(x)\cdot DG=0\\[2mm]
G(x,0)=p\cdot x.
\end{cases}
\ee

\begin{ques}
Does there exist  constant $\overline H(p)$ such that
\be
-\lim_{t\to +\infty}{G(x,t)-p\cdot x\over t}=\overline H(p) \quad \text{uniformly for $x\in \Rset^n$}\ ?  \label{eq:curvature-limit}
\ee
\end{ques}
Due to the presence of curvature term, this question is much harder than the inviscid case ($d=0$).  

First, let us review some literature related to the homogenization of mean curvature type equations.  In other applications such as crystal growth (e.g.  freezing or melting of ice in pure liquid \cite{L1980}), the curvature effect is also considered and  the  motion law is given by (without the drift term $V\cdot \vec n$  and $()_+$ correction)
$$
v_{\vec{n}}=a(x)-d\kappa
$$
for a continuous positive $\Zset^n$-periodic function $a(x)$. The above formula is known as the Gibbs-Thomson relation.
The $()_+$ is not needed in crystal growth since both freezing and melting could occur.  Below is the associated equation
\be
u_t + \left(a(x)- d\, \, \mathrm{div}\left({Du\over |Du|}\right)\right)|Du|=0 \quad \text{in $\Rset^n\times (0,\infty)$}.
\ee

$\bullet$ When $a(x)$ and $Da(x)$ satisfy the following coercivity condition
$$
\min_{\Rset^n}\{a^2-(n-1)d|Da|\}>0,
$$
the existence of effective propagation speed (a corresponding form of \ref{eq:curvature-limit})  and full homogenization have been proved in all dimensions \cite{LS2005} ;

\medskip

$\bullet$  When $n=2$, homogenization was proved in \cite{CM} for all positive $a(x)$ by a geometric approach. See \cite{Caff_13} for a survey on relevant methods. Moreover,  counterexamples are constructed  there when $n\geq 3$. 

\medskip

See \cite{CB, CLS, CN, DKY, GMT, GMOT, Jang, KG} and the references therein for other related works and mathematical models. 
\medskip

The proofs in \cite{LS2005,CM}  rely heavily on the coercivity of HJE in the gradient variable. A new approach needs to be developed to deal with the curvature G-equation where  coercivity is lost when the flow intensity $|V|$ exceeds the local burning velocity $s_L^{0}=1$, which is the typical scenario in turbulent combustion. 

\subsection{Existence of Average Flame Speeds in Cellular Flows}

 Let $V:\Rset^2\to \Rset^2$ be the steady two dimensional cellular flow, e.g, $V(x)=A\,(DH(x))^{\perp}$ with stream function $H(x_1,x_2)=\sin (x_1)\sin(x_2)$ and flow intensity $A>0$. Below is the main theorem proved in \cite{GLXY_22}. 
\begin{thm}\label{cg-thm}
    For any unit vector $p\in \Rset^2$ and initial data $G(x,0)=p\cdot x$, there exists a positive number $\overline H_A(p)$ such that 
\be
\left|G(x,t)-p\cdot x+\overline H_A(p)t\right|\leq C \quad \text{in $\; \Rset^2\times [0,\infty)$}, \label{cG2}
\ee
holds for solution $G(x,t)$ of the curvature G-equation (\ref{cG1}) with a constant $C$ depending only on $d$ and $V$. In particular, this implies that 
$$
-\lim_{t\to \infty}{G(x,t)\over t}=\overline H(p) \quad \text{locally uniformly for $x\in \Rset^2.$}
$$
\end{thm}

The effective Hamiltonian $\overline H_{A}(p)$ corresponds to the effective burning velocity (turbulent flame speed) in the physics literature \cite{Will_85,PR_95,Pet_00}. The inequality (\ref{cG2}) says that to leading order, the solution of curvature G-equation (\ref{cG1}) from the planar initial data $p\cdot x$ develops into a traveling front $p\cdot x-\overline H_A(p)t $.
\medskip

{\bf $\bullet$ Sketch of proof.} Our proof combines PDE methods with a dynamical analysis of the Kohn-Serfaty deterministic game characterization \cite{KS1,KS2} of (\ref{cG1})  
while utilizing the streamline structure of the cellular flow. To the best of our knowledge, this is the first time that a Lagrangian method has been used to prove homogenization of second order PDEs. For the reader's convenience, we will briefly review  key ideas, tools and steps in the proof.  Let us look for a solution to (\ref{cG1}) of the form $p\cdot x-\overline H_A(p) t + v (x)$, where $v$ is the corrector which
satisfies the following equation (a.k.a. cell problem) upon substitution:
\begin{equation}
    \left(1-d\, \; {\rm div}
      \left({p+Dv\over |p+Dv|}\right)\right)_+|p+Dv|+ V(y)\cdot (p+Dv)=\overline H_A(p)
\label{cp1}
\end{equation}
subject to $2\pi$-periodic boundary condition in $y$. The spatial variable of (\ref{cp1}) is renamed $y$ from the original $x$, since  (\ref{cp1}) becomes an independent problem that relates the flow velocity to a large scale constant $\overline H_A$ in the direction $p$. One may view (\ref{cp1}) as a nonlinear eigenvalue problem with $v$ an eigenfunction ($v \pm$ const. still a solution) and $\overline H_A(p)$ an eigenvalue for any given unit vector $p$.
\medskip

As in \cite{LPV}, we start with a  modified (a.k.a. discount) cell problem:
\begin{equation}
  \lambda\, v +  \left(1-d\, \; {\rm div}
      \left({p+Dv\over |p+Dv|}\right)\right)_+|p+Dv|+ V(y)\cdot (p+Dv)= 0,
\label{cp2}
\end{equation}
for a parameter $\lambda >0$, so that the existence and uniqueness of solution $v=v_{\lambda}$ to (\ref{cp2})
is known by Perron's method \cite{CIL_92}. Also comparison principle applies to (\ref{cp2}) where we 
evaluate maximum of $v_\lambda(x)$ to find:
\be\label{v-max}
\max_{x\in \Rset^2}|\lambda \, v_\lambda(x)|\leq 1+\max_{\Rset^2}|V|.
\ee
Our goal is to show that there exists a positive constant $\overline H_A(p)$ such that 
\begin{equation}
\lim_{\lambda\to 0}\lambda \,  v_{\lambda}(x)=-\overline H_A(p) \quad \text{uniformly on $\Rset^2$}, 
    \label{cp3}
\end{equation}
thereby (\ref{cp1}) follows in the 
$\lambda \downarrow 0$ limit with further estimate on $v_{\lambda}$.
\medskip

In the next two subsections, we explain the main ideas for proving (\ref{cp3}) while  omitting the subscript $A$ for simplicity. 
Let $C$ be a constant depending only on $d$ and $V$. 
A {\it key inequality} we shall establish is: 
\begin{equation}
\max_{x,\, y\, \in [-\pi,\pi]^2}|v_\lambda(x)-v_{\lambda}(y)|\leq C,
\label{cp4}
\end{equation}
which implies the constant limit of (\ref{cp3}). 
Due to the curvature term and non-coercivity in the discount cell problem (\ref{cp2}), it is not clear that the equi-continuity of $v_{\lambda}$, a standard yet stronger estimate than (\ref{cp4}), even holds. 
Our strategy is: (1) establishing one-way reachability from the associated game dynamics, (2) 
leveraging a minimum value principle (Eulerian) 
to compensate for the lack of full reachability. The proof coherently integrates Lagrangian and Eulerian thinking.  It turns out that the game trajectory under player I's strategy more or less reverses the propagation route of flame, hence a generalized method of characteristics is at work.

\subsection{Two-Player Game Representation and Analysis}
The deterministic game representation of Kohn-Serfaty \cite{KS1} applies to (\ref{cG1}) on $\Rset^2$ as follows.  Consider the discrete dynamical system $\{x_n\}_{n=1}^{N}\subset\mathbb{R}^2$ associated 
with the game starting from $x_0=x$. For $n=0,1,2,.., N-1$, $|\vec {\eta}_n|\leq 1$ and $b_n\in \{-1,1\}$,
$$
\begin{cases}
x_{n+1}=x_n+\tau\, \sqrt{2d}\, b_n \, \vec{\eta}_n+\tau^2\,\vec{\eta}_{n}^{\perp}-\tau^2\, V(x_n)\\
x_0=x. 
\end{cases}
$$
Player I controls direction via $\vec{\eta}_n$ and player II controls sign via $b_n$. Let $g=g(x)$ be a final payoff function. Player I (II) aims to minimize (maximize) $g(x_N)$. If both players proceed optimally, the value function
$u(x, N\tau^2):=g(x_N),
$
converges to a solution of (\ref{cG1}) with initial data $g(x)$:
$
 \lim_{N\tau^2\to t, \tau\to 0}u(x,N\tau^2)=G(x,t).
$
See Remark \ref{rmk:hgame} for high dimensional version of the game.
\medskip

\begin{rmk}\label{upper-traj}
We shall derive inequalities from {\it one-way reachability in the flow $- V(x)$ via game trajectories} $\{x_n\}_{n=1}^{N}$.
A typical scenario to get an upper bound of the game value is for 
player I to devise a strategy so that the game trajectory, starting at a point $P$, 
ends at a point $Q$ in a desired region $U$ in $N$ moves despite  
any strategy of player II. Then 
\[ u(P, N\tau^2) \leq g(Q) \leq \max_{q\in U} g(q). \]
\end{rmk}
\begin{figure}
   \centering
\includegraphics[scale=0.6]{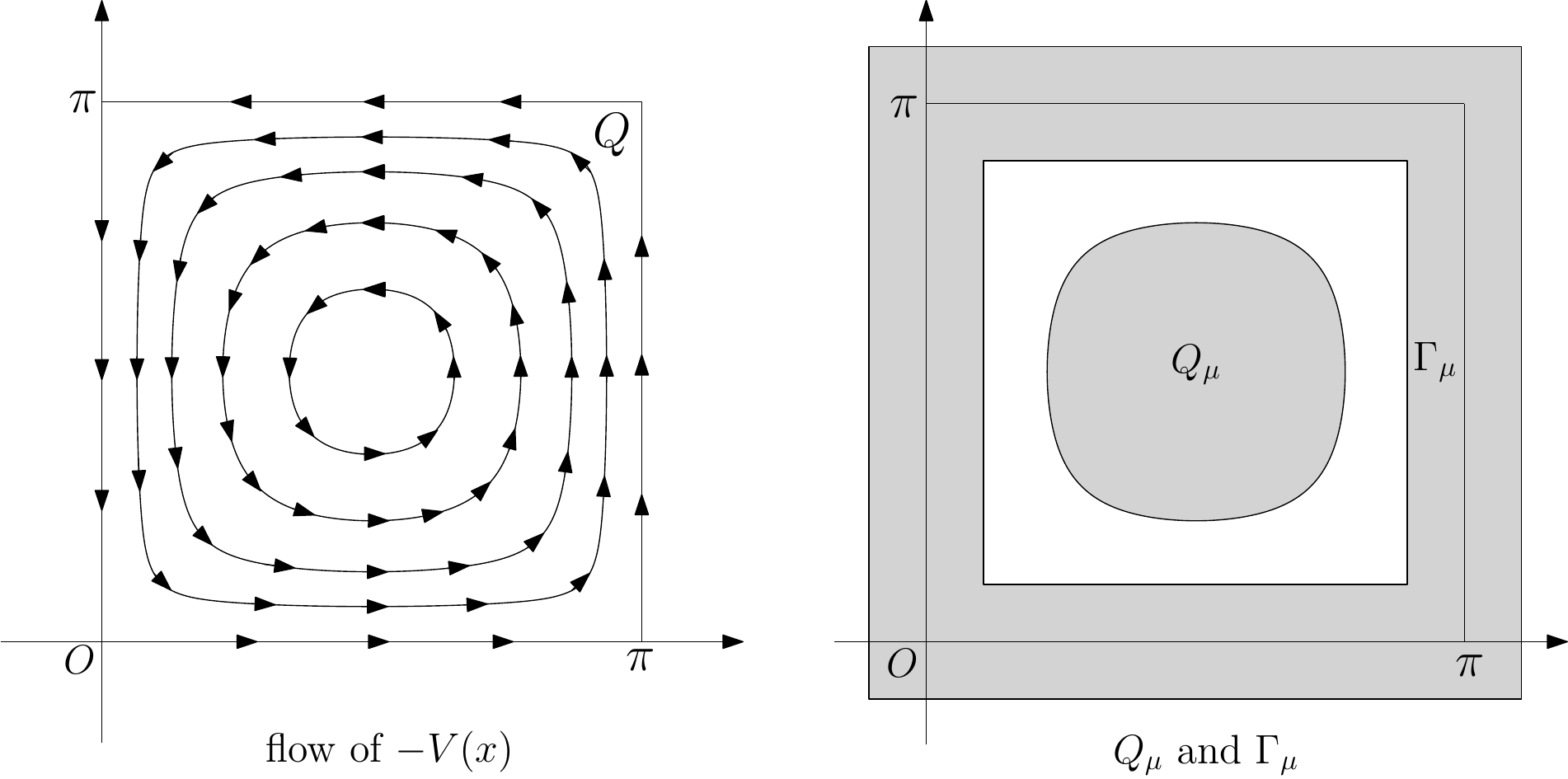} 
\caption{Flow lines of $-V$ (left), two domains in the analysis (right): $Q_{\mu}$ and $\Gamma_{\mu}$}
\label{cquart1}
\end{figure}

\subsubsection{Reachability of Game Trajectory In and Across Cells}
The one-way reachability property of the game dynamics in the interior of cellular flow,  whose quarter-cell with streamline is shown in 
Fig. \ref{cquart1} (left). 
The quarter cell $Q=[0,\pi]\times [0,\pi]$ 
has core 
$Q_\mu=\{x\in Q| \ H(x)>\mu\}$
and boundary region 
$$
\Gamma_\mu=\{x=(x_1,x_2)\in \Rset^2|\ \min\{|x_1|, |x_2|, |x_1-\pi|, |x_2-\pi|\}<\mu\}.
$$
for $\mu\in (0,1]$. Note that $Q_\mu$ is $-V$ flow invariant.

\begin{lem}
Let $\mu \in (0,1)$, then the level curve $\{x\in Q|\ H(x)=\mu\}$  is reachable from $P_1\in Q_\mu$ within time $T_1$ depending only on  $V$.  
Each point $P_2\in \{x\in Q|\ H(x)=\mu\}$ is reachable from $P_1$ within time $C(1+|\log \mu|)$ for a constant $C$ depending only on  $V$.
\label{int-reach}
\end{lem}

\begin{figure}
    \centering
    \includegraphics[scale=0.35]{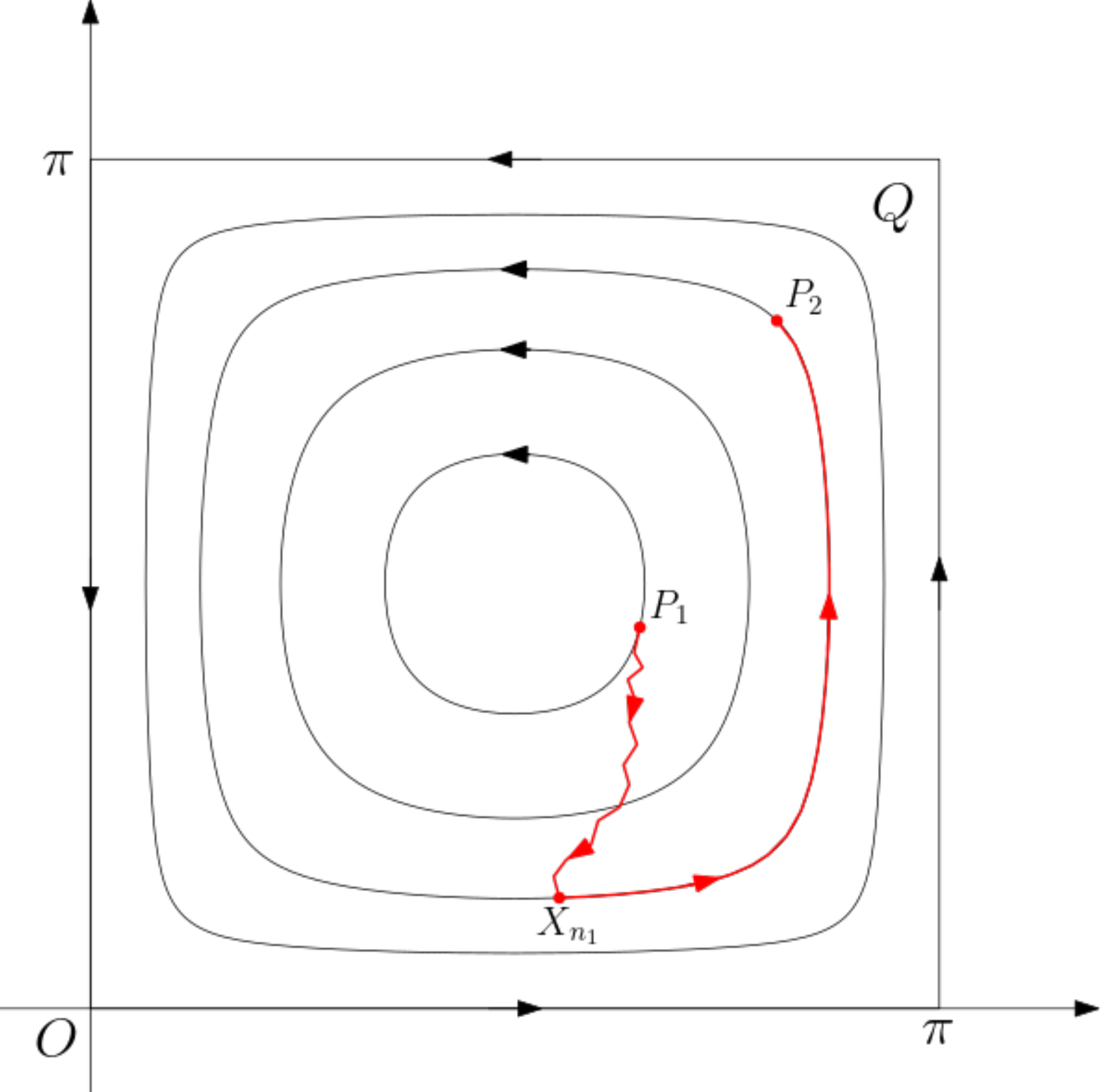} 
    \caption{a game trajectory (in red) of Lemma \ref{int-reach}}
    \label{gametraj1}
\end{figure}


The proof is based on construction of a game trajectory 
in Fig. \ref{gametraj1}.
The reachability 
is only one way. If $P_1$ is near the center of $Q$ where $H(P_1)$ is close to 1 and the curvature on the level curve there exceeds 1, it is NOT reachable from $P_2$. This is different from the inviscid case ($d=0$) where any two points are mutually reachable.

\begin{figure}
    \centering
    \includegraphics[width=0.6\linewidth]{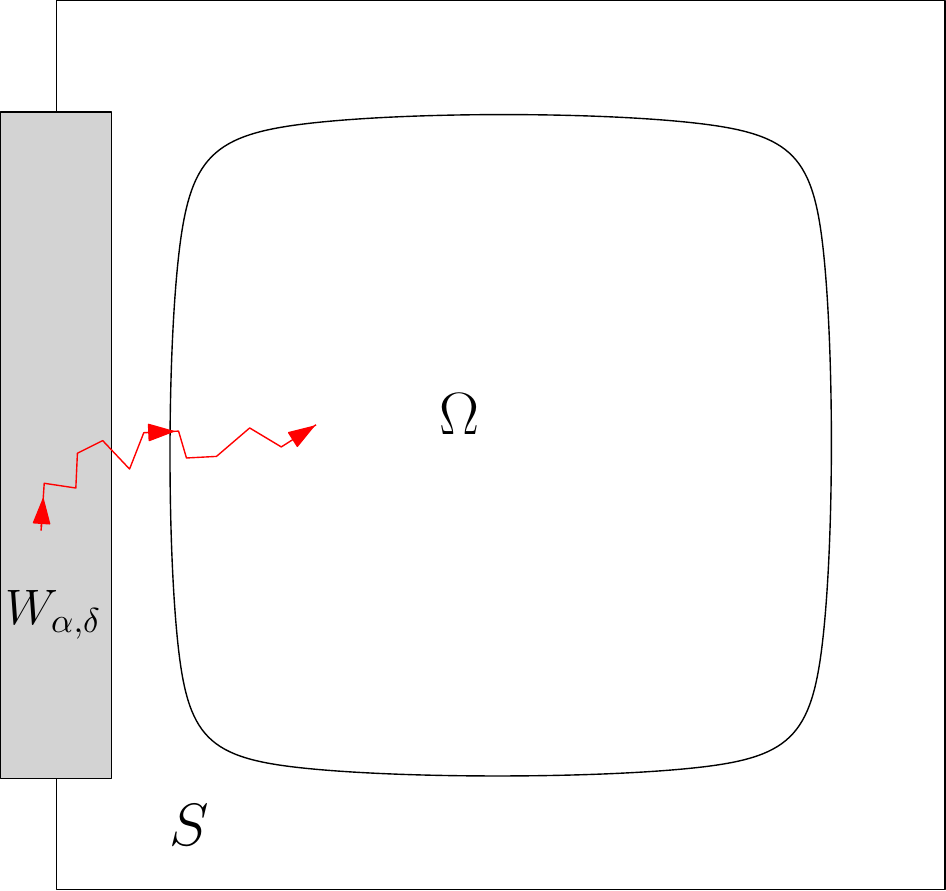}
    \caption{Game trajectory reaching a large enough  interior region $\Omega \subset S=(0,1)^2$ from a thin rectangular boundary layer $W_{\alpha,\delta}=[-\alpha, \alpha]\times [\delta, 1-\delta]$.}
    \label{xcellbdy}
\end{figure}

The game trajectory can go across the quarter-cell boundary and reach any adjacent quarter-cell. More precisely,  
\begin{lem}\label{boundary-reach}
There exists $\mu_0>0$ and $T_{0}>0$ depending only on $d$ and $V$ such that the set $Q_{2\mu_0}$ is reachable from each point $x\in \Gamma_{\mu_0}$ within time $T_{0}$. 
\end{lem}
The crucial part of Lemma \ref{boundary-reach} is that if the 
game trajectory is in the outside portion of the quarter cell boundary layer $\Gamma_{\mu_0}\cap Q^c$, it can pass the boundary and go inside $Q$, as shown in 
Fig. \ref{xcellbdy}, where $S:=(0,1)^2$, $W_{\alpha, \delta}:=[-\alpha, \alpha]\times [\delta, 1-\delta]$, 
for $\alpha>0$ and $\delta\in (0, {1\over 2})$. 
The proof relates an inward optimal game trajectory 
with the outward motion of the zero level set through a PDE comparison argument. 

Based on Lemma \ref{boundary-reach}, one can further show that all quarter cells of cellular flow are reachable by a game trajectory starting from one of them within 
a finite travel time,
which leads to 
\begin{lem}\label{upper-bound} Let $G(x,t)$ be the unique solution of (\ref{cG1}) with $G(x,0)=p\cdot x$. There exist positive constants $\beta$ and $C$ depending only on $d$ and $V$ such that for all $(x,t)\in \Rset^2\times [0, \infty)$,
\be
G(x,t)-p\cdot x\leq -\beta\,  t +C \label{upperbd_G} 
\ee
and 
\be
\max_{x\in \Rset^2}\lambda \, v_{\lambda}(x)< -{\beta \over 2}+ \lambda \, C. \label{upperbd_LdV}
\ee
\end{lem}
The inequality (\ref{upperbd_G}) follows from 
estimating travel times of game trajectories
across cells. The inequality (\ref{upperbd_LdV}) comes from constructing a super-solution to a time-dependent  variant of the discount cell problem (\ref{cp2}) via the inequality 
(\ref{upperbd_G}) where $t$ 
corresponds to $1/\lambda$ 
in (\ref{upperbd_LdV}).
Applying inequality (\ref{upperbd_LdV}), we have from 
Eq. (\ref{cp2}) for $\lambda < \beta/ (4C)$:
\be \label{min-principle}
\left(1-d\, \; \mathrm{div}\left({p+Dv_{\lambda}\over |p+Dv_{\lambda}|}\right)\right)_+|p+Dv_{\lambda}|+ V(y)\cdot (p+Dv_{\lambda})\geq \beta /4
\ee
implying the {\it minimum principle}: the minimum value of $u_{\lambda}:=p\cdot x+v_{\lambda}$ in a domain can only be attained on its boundary. 
\medskip

Next we combine minimum principle with one-way reachability of the game trajectory to prove inequality (\ref{cp4}). First by (\ref{cp2})-(\ref{v-max}), $u= u_{\lambda}$ is a viscosity sub-solution of the stationary G-equation:
\be \label{statgeq}
\left(1-d\, \; \mathrm{div}\left({Du_{\lambda}\over |D u_{\lambda}|}\right)\right)_+|Du_{\lambda}|+ V(y)\cdot (Du_{\lambda}) = 1 +\max_{[-\pi,\pi]^2}\, |V|(y):=\alpha
\ee
for which the following inequality holds:
\be
u_{\lambda}(x_0)\leq  \max_{y\in \bar S}\, u_{\lambda}(y)+\alpha \, T_0 \label{re-ineq}
\ee
if a bounded set $S$ is reachable from $x_0$ via a game trajectory within time $T_0$,  and 
$S$ is invariant under $-V$  flow. To see (\ref{re-ineq}), consider $w := u_{\lambda} - \alpha \, t$ so that $w$ is a sub-solution of G-equation  (\ref{cG1}) with initial data $G(x,0)=u_{\lambda}$. By Remark \ref{upper-traj} and comparison principle, $u_{\lambda}(x_0) - \alpha \, T_0 \leq G(x_0, N\, \tau^2) \leq \max_{y\in \bar S}\, u_{\lambda}(y)$. Note that if the game trajectory reaches $S$ before time $T_0$, then player I can just choose $\vec{\eta}=0$ that makes the trajectory flow within $S$ until $T_0$ due to the flow-invariance of $S$.
\medskip

By Lemma \ref{int-reach} and (\ref{re-ineq}), for each point $x\in \partial Q_{{\mu}}$ and each point $y\in \overline Q_{{\mu}}$, $u_{\lambda}(x)\geq u_{\lambda}(y)-C_{\mu}$ for some constant $C_{\mu}> 0$. Accordingly, 
$$
\min_{x\in \partial Q_{{\mu}}}u_{\lambda}(x)\geq \max_{x\in \overline Q_{{\mu}}}u_{\lambda}(x)-C_{\mu}.
$$
By minimum principle:
$
\min_{x\in \partial Q_{{\mu}}}u_{\lambda}(x)=\min_{x\in \overline Q_{{\mu}}}u_{\lambda}(x)
$, and so:
\be\label{interior-control}
\max_{x\in \overline Q_{{\mu}}}u_{\lambda}(x)-\min_{x\in \overline Q_{{\mu}}}u_{\lambda}(x)=\max_{x,y\in \overline Q_{{\mu}}} |u_\lambda(x)-u_\lambda(y)|\leq C_{\mu}.
\ee
With the help of minimum principle and more delicate analysis, similar estimates hold over the cell boundary and in other quarter cells, and so the key inequality (\ref{cp4}) holds. Theorem \ref{cg-thm} follows with a sub/super solution argument \cite{GLXY_22}. 
\medskip

\subsection{Existence of $\overline H$ for two dimensional (2D) incompressible flows}

    

The proof for the case of cellular flow in the last subsection can be modified to establish the existence of $\overline H(p)$ for a large class of 2D incompressible flows \cite{GLXY_22}, for example, if $V=D^{\perp}H$ for a smooth periodic stream function $H$ that has finitely many non-degenerate critical points whose flow structure is well described in \cite{A1991}. Also, the constant local burning velocity can be replaced by a positive, continuous and periodic function $a(x)$. 

\medskip

\begin{ques}
Does  (\ref{eq:curvature-limit}) hold for all 2D periodic  incompressible flows?
\end{ques}

\subsection{Bifurcation of 
effective burning velocity in 3D shear flows}

A natural question is whether our robust approach for 2D incompressible flows can be extended to 3D flows.  Surprisingly, the effective burning velocity ceases to exist when the flow intensity surpasses a threshold value (bifurcation) for general 3D shear flows, the simplest class of 3D incompressible flows.  Meanwhile, the existence of effective burning velocity for 2D shear flows can be easily established  by PDE methods.  The methods in this section mainly rely on PDE approaches. It is not clear to us how to interpret it by 3D  game theory (see Remark \ref{rmk:hgame}).

Assume that $f$ is a Lipschitz continuous function satisfying 
\be\label{c1}
 \{x\in \Rset^2\,:\, f(x)=\max_{\Rset^2}f\}=\Zset^2, \;  \{x\in \Rset^2\,:\, f(x)=\min_{\Rset^2}f\}=\Zset^2+\vec{q}
\ee
for some $\vec{q}\in \Rset^n$.  Here we do not pursue the optimal assumptions on  $f$.   Consider the 3D shear flow
$$
V(x)=(0, 0, Af(x')),
$$
for $x=(x',x_3)\in \Rset^3$ and $x'\in \Rset^2$, the constant $A>0$ represents the flow intensity. Then equation (\ref{eq:linearcur}) is reduced to 
\be
\begin{cases}
v_t+\left(1-d\,{\rm div}{\frac{p+Dv}{\sqrt{p_{3}^2+|p'+Dv|^2}}}\right)_+\sqrt{p_{3}^2+|p'+Dv|^2}+A p_{3} f(x')=0,\\[3mm]
v(x',0)=0
\end{cases}
\ee
for $G(x,t)=v(x',t)+p\cdot x$ and $p=(p',p_3)\in \Rset^3$ and $p'\in \Rset^2$.

For any $p\in \Rset^3$ with $p_3\not=0$,  denote by $S_H$ the set of all $A\geq 0$ such that 
\be\label{eq:Llimit}
\lim_{t\to \infty}-{G(x,t)\over t}=\overline H(p,A) \quad \text{for all $x\in \Rset^3$}. 
\ee
for a constant $\overline H(p,A)$. 

For $\lambda>0$, let $v_\lambda$ be the periodic continuous viscosity solution to the discount $\lambda$-problem on $\Rset^2$:
$$
\lambda v_{\lambda}+ \left(1-d\,{\rm div}{\frac{p+Dv_\lambda}{\sqrt{p_{3}^2+|p'+Dv_\lambda|^2}}}\right)\sqrt{p_{3}^2+|p'+Dv_\lambda|^2}+A p_{3} f(x')=0.
$$
Then (\ref{eq:Llimit}) is equivalent to 
$$
\lim_{\lambda\to 0}\lambda \, v_{\lambda}(x')=\overline H(p,A), \quad \text{for $x'\in \Rset^2$}.
$$
The following result is proved in \cite{MMTXY2023}:
\begin{thm} For any $p\in \Rset^3$ with $p_3\not=0$, 
$$
S_H(p)=[0,A_0(p)]
$$
for a constant $A_0(p)$ that is  continuous with respect to $p$ and $A_0(sp)=A_0(p)$ for any $s>0$. 
\end{thm} 
   
Consequently,   $\bar A=\min_{p\in \Rset^3} A_0(p)$ is the bifurcation value for the full homogenization of the corresponding curvature G-equation with uniformly continuous initial data. 

Below, we explain how to characterize the bifurcation value $A_0(p)$.  Consider the cell problem without $()_+$ cutoff:
\be\label{eq:cellnocut}
\left(1-d\,{\rm div}{\frac{p+D\tilde v}{\sqrt{p_{3}^2+|p'+D\tilde v|^2}}}\right)\sqrt{p_{3}^2+|p'+D\tilde v|^2}+A p_{3} f(x')=\overline {\tilde H}(p,A). 
\ee
By the Bernstein technique,  it can be shown that, unlike the physical cut-off case, for any given $p\in \Rset^3$,  there exists a constant $\overline {\tilde H}(p,A)$ such that  the above equation always has a $C^{2,\alpha}$ periodic solution $\tilde v=\tilde v(x)$. This can be viewed as a special case of \cite{LS2005}. Moreover, $\overline {\tilde H}(p,A)-F(p)A$ is strictly decreasing  for $A\geq 0$ and $F(p)=\max_{x'\in\Rset^2}(p_3f(x'))$. 

Apparently, if  $\overline {\tilde H}(p,A)\geq AF(p)$, the $\tilde v$ is also a solution to the equation with physical cut-off: 
$$
\left(1-d\,{\rm div}{\frac{p+D\tilde v}{\sqrt{p_{3}^2+|p'+D\tilde v|^2}}}\right)_+\sqrt{p_{3}^2+|p'+D\tilde v|^2}+A p_{3} f(x')=\overline {\tilde H}(p,A). 
$$
Then $A\in S_H(p)$ and 
$$
\lim_{t\to \infty}-{G(x,t)\over t}=\overline {\tilde H}(p,A).
$$

{\bf Characterization of $A_0(p)$}. For fixed $p\in \Rset^3$, the bifurcation value $A_0(p)$ is the unique number such that
\be\label{eq:chara}
\overline {\tilde H}(p,A_0(p))=F(p)A_0(p).
\ee
Moreover,  
$$
\overline { H}(p,A)=\overline {\tilde H}(p,A) \quad \text{for $A\in [0,A_0(p)]$}.
$$

Hereafter, we reviewed main ideas of the proof.  For convenience, we assume that $\max_{\Rset^2}f=0$ and $p_3=1$. Then $F(p)=0$. 

\medskip

{\bf Step 1:} Using comparison principle,  we can deduce that the set $s_H$ is either  $[0,\infty)$ or $[0,A_0(p)]$ for a positive constant $A_0(p)$.

\medskip

{\bf Step 2:} We need to show that when $A$ is large enough, $\overline {\tilde H}(p,A)<0$. The main idea is that if that is not true, then we can construct a radial symmetric  function $v(x)=g(|x|)$ satisfying 
$$
\left(1-d\,{\rm div}{\frac{p+Dv}{\sqrt{1+|p'+Dv|^2}}}\right)\sqrt{1+|p'+Dv|^2}\geq 1 \quad \text{for $r\leq |x|\leq 2r$}
$$
for small $r$ or equivalently in terms of $g:[r,{1\over 2}]\to \Rset$
$$
\sqrt{1+(g')^2}-d\left({g''\over 1+(g')^2}+{g'\over r}\right)\geq 1 \quad \text{for $r\in \left(r_0,\  3r_0\right)$}.
$$
Both are in the viscosity sense. Moreover, $g(2r)-g(r)\geq 1$. This will lead to a contradiction when $r$ is small due to the presence of the term ${g'\over r}$ on the left hand side. We note in passing that for the 2D shear flow, such a term is absent. 

\medskip

{\bf Step 3:}  Finally we need to verify the characterization (\ref{eq:chara}),  which will provide the continuous dependence on $p$. It suffices to show that 
if $\lim_{\lambda\to 0}\lambda v_{\lambda}(x)=0$, then $\overline {\tilde H}(p,A)=0$.  This can be done in two steps. Let $u_{\lambda}(x)=v_{\lambda}(x)+p'\cdot x$.

\medskip

 {\bf Step 3.1} Employing estimates of minimal surface type equations established in \cite{S1}, we can show that as $\lambda\to 0$,  $u_{\lambda}-u_{\lambda}(0)$ in $\Rset^2\backslash\Zset^2$ locally uniformly converges to a function $u$ that is a $C^{2,\alpha}$ solution to 
 $$
 \left(1-d\,{\rm div}{\frac{Du}{\sqrt{1+|Du|^2}}}\right)\sqrt{1+|Du|^2}+Af(x')=0 \quad \text{in $\Rset^2\backslash\Zset^2$}. 
 $$

 \medskip

 {\bf Step 3.2} Using integration by parts through suitable test functions, the isolated singularities $\Zset^2$ are proved to be removable at least in the viscosity sense. As a result,  comparing with the cell problem without cutoff (\ref{eq:cellnocut}) leads to $\overline {\tilde H}(p,A)=0$.

 \medskip

 \begin{rmk}\label{rmk:hgame} The Kohn-Serfaty game can be naturally extended to higher dimensions \cite{KS1}. Consider the three dimensional (3D) case. For $k\geq 1$, at the $k$-th step,  let $\tau$ be the time step size. 
\medskip

1. Player I  chooses three  vectors $u_k$,  $v_k$ and $w_k$ such that  they are mutually perpendicular and $|u_k|=|v_k|=|w_k|\leq 1$.
\medskip

2. Player II  chooses $b_k=\pm 1$ and $c_k=\pm 1$ to replace ($v_k$,$w_k$) by ($b_kv_k$, $c_kw_k$).
\medskip

3. Update: $X_{k}=X_{k-1}+\tau \sqrt{2d}\, (b_kv_k+c_kw_k)+\tau^2u_{k}-\tau^2V(X_{k-1})$.

\medskip

The above generalization is based on the following fact: given a symmetric matrix $S$ and mutually orthogonal unit vectors $u$, $v$ and $q$, 
$$
u\cdot S\cdot u+v\cdot S\cdot v+q\cdot S\cdot q=tr(S).
$$
The other parts are similar to the 2D case.  A major difficulty of the 3D  game is that it is hard for play I to design a strategy to move between orbits of $\dot \xi=-V(\xi)$ due to the lack of topological restrictions. It is not clear to us if the bifurcation phenomenon in shear flows is typical for 3D flows.
     
 \end{rmk}

The following question is the  first step towards understanding the role of game theory in 3D flows.
\medskip

\begin{ques}
Can we find a game theoretic interpretation of the failure of homogenization in 3D shear flows? 
\end{ques}
\medskip
Unlike shear flows,  turbulent 3D flows often have swirl structures (eddies), which help mix things well.  Accordingly, an interesting and challenging problem is: 

 \medskip

 \begin{ques}
 If $V$ is the ABC flow,  do we have 
 $$
S_H(p)=[0,\infty)
$$
for all $p\in \Rset^3$?
\end{ques}
\subsection{Impact of curvature on effective burning velocity}

A significant problem in the combustion literature is to understand how the curvature term impacts the prediction of the burning velocity. There is a consensus that curvature slows down flame propagation \cite{PR_95}.  Heuristically, this is because the curvature term smooths out the flame front and reduces the total area of chemical reaction. Under the G-equation model, this is equivalent to showing that the effective burning velocity $\overline H$ is decreasing with respect to the Markstein number $d$, which has been observed in  combustion  experiments \cite{CWL2013} and numerical computations \cite{KAW1988, EnsKin_23}. We would like to point out that different Markstein number is achieved by mixture of different fuels in experiments \cite{CWL2013}.  To establish it rigorously from curvature G-equation is challenging.  Below we look at shear flows and 2D cellular flows.

\subsubsection{Decrease with respect to the Markstein number in  shear flows} The following is the first mathematically rigorous result in this direction obtained in \cite{LXY_18} for two dimensional shear flows.

\begin{thm} Let $V(x)=(0,f(x'))$ for $x=(x_1,x')\in \Rset^2$ and $f$ is a non-constant periodic continuous function. Then for any fixed vector $p=(p_1,p_2)\in \Rset^2$ with $p_2\not=0$ and $h(d)=\overline H(p,d)$,
$$
h'(d)\leq 0  \quad \text{for\;  $d>0$}
$$
and
$$
h'(d)< 0  \quad \text{when $d$ is close to 0}.
$$
In particular, this implies that
$$
\overline H(P,d)<\overline H(P,0).
$$
\end{thm}

\noindent {\bf Sketch of proof.} For 2D shear flows, the effective burning velocity $\overline H(p,d)$ always exists and 
$$
\overline H(p,d)=\max\left\{\overline {\tilde H}(p,d),\ F(p)\right\}.
$$
If $d$ is small, $\overline H(p,d)=\overline {\tilde H}(p,d)
$. Here $F(p)=\max_{x'\in\Rset}(p_2f(x'))$ and $\overline {\tilde H}(p,d)$ is from (\ref{eq:cellnocut}) without $()_+$ correction.  Hence it is enough to show that  for fixed unit vector $p$, $\overline {\tilde H}(p,d)$ is strictly decreasing with respect to $d$

\medskip

{\bf Step 1} For convenience, write
$$
E(d)=\overline {\tilde H}(p,d). 
$$
Without loss of generality, we may assume that $p_2=1$. Then the cell problem is reduced to an ODE satisfied by $\phi(x')=p_1+\tilde v'(x')$ for $x'\in \Rset$
$$
-{d\phi'\over 1+\phi^2}+\sqrt{1+\phi^2}+f(x')=E(d). 
$$
Taking derivative on both sides with respect $d$, we obtained a linear ODE satisfied by $F=\phi_d$,
$$
-dF'(x')+b(x')F=E'(d)(1+\phi^2)+\phi'
$$
Here $b(x')={2d\phi'\phi\over 1+\phi^2}+\phi\sqrt{1+\phi^2}$. 

\medskip

{\bf Step 2}: Due to the periodicity of $F$,  $E'(d)$ can be expressed by complicated integrations involving $\phi$ and $\phi'$.  Then $E'(d)<0$ can be  proved by establishing the following delicate inequality
$$
\begin{array}{ll}
&e^{T}\int_{0}^{T}f(t)\, e^{-t}\int_{0}^{t}g(f(s))\, e^s\,ds\, dt+\int_{0}^{T}f(t)\, e^{-t}\int_{t}^{T}g(f(s))\, e^s\,ds\, dt\\[5mm]
&\geq (e^T-1)\int_{0}^{T}f(t)g(f(t)))\,dt+{\theta\over 2} \int_{[0,T]^2}|f(t)-f(s)|^2\,dt\, ds.
\end{array}
$$
for any continuous positive function $f\in C([0,T])$ and $g\in C^1((0,L])$ for $L=\max_{[0,T]}f$ satisfying $g'\leq -\theta$ for some $\theta\geq 0$. 

\medskip

{\bf Step 3:} To prove the above inequality, we look at a suitable discrete form where the question is reduced to finding maximum points of functions of finitely many variables. Hence the conclusion can be derived by the simple fact that the gradient vanishes at maximum points.

\medskip

Note that the proofs in \cite{LXY_18} heavily depend on 1D (or $x'\in \Rset$).  Hence a natural question is 

\medskip

\begin{ques}
Consider 3D shear flow $V(x_1,x_2,x_3)=(0.0,f(x_1,x_2))$ and the associated $E(d)=\overline {\tilde H}(p,d)$ for $p=(p',1)\in \Rset^3$.  Do we have  
$$
E'(d)<0
$$
for $d>0$ and non-constant $f$?
\end{ques}
\subsubsection{Slow down in two dimensional cellular flows} 
\begin{figure}
\begin{center}
\includegraphics[scale=0.6]{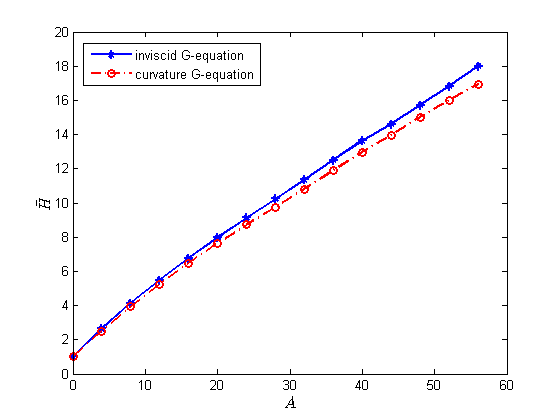}
\caption{Comparison of the homogenized Hamiltonian function in the cellular flow amplitude $A$ without curvature ($\overline H(p,0,A)$, blue line-star) and with curvature ($\overline H(p,0.1, A)$, red dashed-line-circle) by a finite difference method \cite{LXY_13a}.}
\label{average-front}
\end{center}
\end{figure}

In case of 2D cellular flow, numerical computation at $d=0.1$ and different values of $A$ (Fig. \ref{average-front})  suggests that for fixed flow intensity 
\be\label{eq:onewaybound}
\overline H(p,d,A)<\overline H(p,0,A)
\ee
and the gap
$$
\overline H(p,0,A)-\overline H(p,d,A)
$$
increases as $A$ increases. 

It remains a very interesting question to prove (\ref{eq:onewaybound})  rigorously. Due to the lack of smooth solutions,  we cannot take the derivative of the equation to track the change with respect to $d$ as for the shear flows,  i.e., the PDE method in the previous section does not work.  Accordingly,  the only available tool we have  is the game theory interpretation.  However, to prove (\ref{eq:onewaybound}) might require  accurate analysis of the game trajectory, which is very difficult, in particular,  given the hidden nonlinear stochastic nature of the game.  One hopes that when the flow intensity $A$ gets large (the typical scenario in turbulent combustion), the subtle curvature effect might be amplified and become easier to capture. Thus,  a more realistic  question is 

\medskip

\begin{ques}
Do we have 
$$
\lim\inf_{A\to \infty}(\overline H(p,0,A)-\overline H(p,d,A))>0\ ?
$$
Or do we actually have the limit equal to  $\infty$? 
\end{ques}
We can show that $\overline H(p,d,A)$ also grows like $O(A/\log A)$ as in the inviscid case ($d=0$). Hence, finding the constant (if it exists) in Question 2 seems relevant.  
\medskip

\section{Strain G-equation}
\setcounter{equation}{0}
Recall that G-equation with a strain term has the following form:
\be\label{strainG}
G_t+\left(1+d\ {DG\cdot S(x)\cdot DG\over |DG|^2}\right)_+|DG|+V(x)\cdot DG=0,
\ee
where $S={DV+(DV)^{\top}\over 2}$ is the strain rate tensor.  This is a non-coercive and non-convex Hamilton-Jacobi equation (HJE). We again consider the existence of effective burning velocity and its dependence on the constant $d$ and flow intensity. It is well-known that solutions to non-convex first-order Hamilton-Jacobi equations can be interpreted by deterministic two-person differential games \cite{I1965}, which will play important roles in our analysis.  

\subsection{Shear flows}   Consider the $n+1$-dimensional shear flow
$$
V(x)=(0,0,...,0, f(x'))
$$
for $x=(x',x_{n+1})\in \Rset^{n+1}$ and $x'\in \Rset^n$. Then for $p=(p',1)\in \Rset^{n+1}$,  the  corresponding cell problem is reduced to 
$$
\left(1+{d(p'+Dv)\cdot Df\over {1+|p'+Dv|^2}}\right)_+\sqrt{{1+|p'+Dv|^2}}+f(x')=\overline H(p',d).
$$
First,  let $\overline {\tilde H}(p',d)$ be the effective Hamiltonian when the $()_+$ correction is absent: 
$$
\left(1+{d(p'+D\tilde v)\cdot Df\over {1+|p'+D\tilde v|^2}}\right)\sqrt{{1+|p'+D\tilde v|^2}}+f(x')=\overline {\tilde H}(p',d).
$$
The existence follows directly from \cite{LPV}. Since $Df=0$ when $f$ attains maximum,  it is easy to see that for all $d\geq 0$
$$
\min_{p'\in \Rset^n}\overline {\tilde H}(p',d)\geq 1+\max_{\Rset^n}f.
$$
 Hence $\tilde v$ is also a solution to 
$$
\left(1+{d(p'+D\tilde v)\cdot Df\over {1+|p'+D\tilde v|^2}}\right)_+\sqrt{{1+|p'+D\tilde v|^2}}+f(x')=\overline {\tilde H}(p',d).
$$
Accordingly, $\overline {\tilde H}(p',d)=\overline {H}(p',d)$, i.e., the cell problems with or without the $()_+$ corrections are the same for shear  flows. Hence we write 
\be\label{eq:strainshear}
\left(1+{d(p'+Dv)\cdot Df\over {1+|p'+Dv|^2}}\right)\sqrt{{1+|p'+Dv|^2}}+f(x')=\overline H(p',d).
\ee
Here the corresponding Hamiltonian is 
$$
H(P,y)=\sqrt{1+|P|^2}+{dP\cdot Df(y)\over \sqrt{1+|P|^2}}
$$
for $(P,y)\in \Rset^n\times \Rset^n$. Note that $H$ is not convex in the $P$ variable.  When $n=1$, $H$ is quasiconvex is $P$, which however is not true when $n\geq 2$.

An interesting question is how the presence of the strain rate term impacts the value of effect burning velocity $\overline H$. Numerical results show that $\overline H$ is decreasing with respect to the Markstein number $d$ \cite{KAW1988}. This is consistent with the following theorem that says that for two dimensional shear flows, it reduces the value of $\overline H(p',d)$

\begin{thm}\label{thm:strainslow} Assume that $n=1$ and fix $p'\in \Rset$. Then $\overline H(p',d)$ is strictly decreasing with respect to $d\geq 0$ when $\overline H(p',d)>1+\max_{\Rset}f=\min_{p\in \Rset^n}\overline {H}(p,d)$.
\end{thm}

Proof: For simplicity, write $E(d)=\overline {H}(p',d)$ and $u(x)=p'x+v(x)$ for $x\in \Rset$. Then
\be\label{eq:1dstrain}
\left(1+{du'f'\over {1+|u'|^2}}\right)\sqrt{{1+|u'|^2}}+f(x')=E(d).
\ee
Suppose that $E(d)>1+\max_{\Rset}f$.  Note that  for fixed $x\in \Rset$, the Hamiltonian is 
$$
H(p,x)=\sqrt{1+p^2}+{dpf'(x)\over \sqrt{1+p^2}}
$$
for $(p,x)\in \Rset\times \Rset$. Then  the derivative of $H$ with respect to $p$ has the form:
$$
H_p(p,x)={p\over \sqrt{1+p^2}}+{df'(x)\over (\sqrt{1+p^2})^3}={p(1+p^2)+df'(x)\over (\sqrt{1+p^2})^3}.
$$
Apparently, for fixed $x$, there exists a unique $p=p(x)$ such that $H_p(p(x),x)=0$. Then $H$ attains minimum at $p=p(x)$ and $H_p(p, x)<0$ when $p<p(x)$ and $H_p(p, x)>0$ when $p>p(x)$.
Thus $H^{-1}$ is a well defined smooth function on either $(-\infty, p(x)]$ or $[p(x),\infty)$. Since 
$$
\min_{p\in \Rset}H(p,x)\leq H(0,x)=1< E(d)-f(x),
$$
we have either 
$$
u'(x)>\max\{0,p(x)\} \quad \text{for a.e. $x\in \Rset$}
$$
or
$$
u'(x)<\min\{0,p(x)\} \quad \text{for a.e. $x\in \Rset$}.
$$
Therefore $u$ must be smooth due to the strict monotonicity of $H$ on either side of $p=p(x)$. Without loss of generality, we assume that 
$$
u'(x)>\max\{0,p(x)\} \quad \text{for all $x\in \Rset$}.
$$
Then for fixed $x$
$$
H_p(u'(x),x)>0 \quad \mathrm{and}\quad u'(x)=H^{-1}(E(d)-f(x),x).
$$
Here for fixed $x$, $H^{-1}(\cdot,x)$ is the inverse map of $H(\cdot,x)$ for $p\geq p(x)$. Let $g(x)=u'(x)>0$. Taking derivative of (\ref{eq:1dstrain}) with respect to $d$ says that
$$
H_p(g(x),x)g_d+{gf'\over \sqrt{1+g^2}}=E'(d).
$$
Since $h(x)={1\over g(1+g^2)}>0$ and
$$
H_p(g(x),x)={g\over \sqrt{1+g^2}}(1+dh(x)f'(x))>0,
$$
we have that
$$
g_d={E'(d)-{gf'\over \sqrt{1+g^2}}\over H_p}.
$$
Thanks to $\int_{0}^{1}g_d(x)\,dx=0$, 
$$
E'(d)\int_{0}^{1}{1\over H_p(g(x),x)}\,dx=\int_{0}^{1}\frac{{gf'\over \sqrt{1+g^2}}}{H_p}\,dx.
$$
Then
$$
\frac{{gf'\over \sqrt{1+g^2}}}{H_p}={f'(x)\over 1+dh(x)f'(x)}
$$
So
\begin{align*}
\int_{0}^{1}\frac{{gf'\over \sqrt{1+g^2}}}{H_p}\,dx&=\int_{0}^{1}{f'(x)\over 1+dh(x)f'(x)}\,dx\\
&=\int_{0}^{1}{f'(x)\over 1+dh(x)f'(x)}\,dx-\int_0^1f'(x)\,dx\\
&=-\int_{0}^{1}{dh(x)(f'(x))^2\over 1+dh(x)f'(x)}\,dx<0.
\end{align*}
\qed

\medskip

\begin{ques}
Does Theorem \ref{thm:strainslow} hold for $n=2$ (or 3D shear flows) ?
\end{ques}

\subsection{2D cellular flows} Let
\be\label{vortex}
V(x)=(-H_{x_2},H_{x_1})
\ee
for $x=(x_1,x_2)\in \Rset ^2$ and $H=\sin x_1\, \sin x_2$. The corresponding strain tensor is
$S=
\left (\begin{array}{ccc}
-\Phi& 0\\
0& \Phi
\end{array}\right)
$ for $\Phi(x)=\cos x_1\, \cos x_2$.  The strain G-equation with linear initial data then takes form of
\be\label{planar}
\begin{cases}
G_t+\left(1-Ad \Phi(x) {|G_{x_1}|^2-|G_{x_2}|^2\over |DG|^2}\right)_{+}|DG|+AV(x)\cdot DG=0\\
G(x,0)=p\cdot x.
\end{cases}
\ee
Here $A\geq 0$ represents the flow intensity. We expect the existence of effective burning velocity for cellular flows to be established by a strategy similar to that in section 4 using two-player differential game (revisited in subsection 5.3): (1) first prove some one way reachability, and then (2) use the minimum value principle to take care of the other direction.  The argument will be simpler for the strain G-equation due to the deterministic nature of its game dynamics valid in all dimensions. On the other hand, compared with the curvature G-equation,  one analytical disadvantage of the strain G-equation is the lack of clear connection between the value of the strain term and the geometric structure of the streamlines. 

\medskip
\begin{ques}
Does the effective burning velocity exist in (\ref{planar})  for general incompressible periodic flows in all dimensions?
\end{ques}

\medskip

Expressing  the associated Hamiltonian 
$$
H(p,x)=\left(1-Ad \Phi(x) {|p_1|^2-|p_2|^2\over |p|^2}\right)_{+}|p|+AV(x)\cdot p
$$
in a min-max or max-min form of differential game (see (\ref{eq:player1}) or (\ref{eq:player2}) in subsection 5.3) is not easy. It is helpful to first consider a simplified Hamiltonian such as
$$
\tilde H(p,x)=\left(|p|-Ad \Phi(x) (|p_1|-|p_2|)\right)_{+}+AV(x)\cdot p
$$
when applying a game theory interpretation and  then employ comparison principle to control solutions to the strain G-equation. 

For mathematical interest,   let us also consider the equation without $``+"$ correction
$$
\begin{cases}
\tilde G_t+\left(1-Ad \Phi(x) {|\tilde G_{x_1}|^2-|\tilde G_{x_2}|^2\over |D\tilde G|^2}\right)|D\tilde G|+AV(x)\cdot D\tilde G=0\\
\tilde G(x,0)=p\cdot x.
\end{cases}
$$
 Then it was proved in \cite{XY_14a} that there are different ranges of $A$ in terms of the existence and non-existence of effective Hamiltonian. This is quite different from the strain G-equation where the existence of effective burning velocity is expected to hold for all $A>0$ . For simplicity,  we state the result for $p=e_1=(1,0)$ below.
\begin{thm}  There exists $0<A_0<A_1$ such that

(i) (Homogenization and propagation range) If $A\in [0,A_0)$ for $A_0={1\over d\max_{x\in \Rset^n}||S(x)||}$, then the effective Hamiltonian
$$
\lim_{t\to \infty}-{\tilde G(x,t)\over t}=\overline H(e_1,A)  \quad \text{for all $x\in \Rset^2$}
$$
exists and is positive. The proof is similar to that of the basic case ($s_l=1$) \cite{XY_10}. 

\medskip

(ii)  (Homogenization and trapping range) There exists $A_1>0$ such that when $A\geq A_1$, 
$$
\lim_{t\to \infty}-{\tilde G(x,t)\over t}=0 \quad \text{for all $x\in \Rset^2$}.
$$
The proof is based on a very delicate analysis of the game dynamics (in subsection 5.3) and the streamline structure of cellular flow. 

\medskip

(iii) (Intermediate non-homogenization range) There exists $A\in [A_0, A_1)$ such that the effective Hamiltonian does not exist, i.e., 
$$
\text{ $\lim_{t\to \infty}-{\tilde G(x,t)\over t}$, if it exists, depends on $x$}. 
$$
More precisely,  the limit (if it exists) is non-positive when $x$ is very close to integer points $\Zset^2$ and is positive when $x$ is away from those integers points. 
\end{thm}

\subsection{Two-person differential game for 1st order nonconvex HJEs} We revisit the two person zero sum 
differential game representation for first order non-convex Hamilton-Jacobi equations (HJEs) \cite{I1965} based on  
notations and formulations from \cite{ES1984}.  For simplicity, we only consider HJEs for front propagation:
$$
\begin{cases}
u_t+H(Du,x)=0,  \quad \text{on $\Rset^n\times (0,\infty)$},\\
u(x,0)=g(x),
\end{cases}
$$
and refer to \cite{ES1984} for general HJEs. 
The equations in \cite{ES1984} are formulated as terminal value problems. To be consistent with other parts of the paper,  we adopt the initial value problem.  Assume that
\be\label{eq:player1}
H(-p,x)=\max_{y\in Y}\min_{z\in Z}\{f(x,y,z)\cdot p\} \quad \text{for all $(x,p)\in \Rset^n\times \Rset^n$}
\ee
or
\be\label{eq:player2}
H(-p,x)=\min_{z\in Z}\max_{y\in Y}\{f(x,\eta,\mu)\cdot p\} \quad \text{for all $(x,p)\in \Rset^n\times \Rset^n$}.
\ee
Here $Y$ and $Z$ are subsets of $\Rset^n$. \nit We also set

(1) $M(t)$ as the set of measurable functions $[0,t]\to Y$;

(2) $N(t)$ as the set of measurable functions $[0,t]\to Z$;

(3) $\Gamma (t)$ as the set of strategies of  player I, i.e, non-anticipating mappings $\alpha: N(t) \to M(t)$ satisfying that for all $s<t$ and $z,\tilde z\in N(t)$:
$$
\begin{cases}
z(s)={\tilde z}(s)  \quad \text{for a.e.  $0\leq s\leq t$}\\
\mathrm{implies\ that}\ \alpha (z)(s)=\alpha({\tilde z})(s) \quad \text{for a.e.  $0\leq s\leq t$};
\end{cases}
$$
(4) $\Delta (t)$ as the set of strategies of  player II, i.e., non-anticipating mappings  $\beta: M(t)\to N(t)$ satisfying that for all $s<t$ and $y,\tilde y\in M(t)$:
$$
\begin{cases}
y(s)={\tilde y}(s)  \quad \text{for a.e.  $0\leq s\leq t$}\\
\mathrm{implies\ that}\ \ \beta(y)(s)=\beta({\tilde y})(s) \quad \text{for a.e.  $0\leq s\leq t$}.
\end{cases}
$$
Given $t>0$ and $(y(s), \beta)\in M(t)\times \alpha(t)$ or $(z(s),\alpha) \in N(t)\times \beta(t)$, let $x=x(s):[0,t]\to \Rset ^2$ be a Lipschitz continuous curve satisfying $x(0)=x$ and

(i) if at each time $s$,  player I moves first by choosing $y(s)$  and player II chooses a corresponding strategy $\beta$, then
$$
\dot x(s)=f(x(s),  y(s),  \beta(y)(s)) \quad \text{for a.e $s\in (0,t)$};
$$

(ii) if at each time $s$, player II moves first by choosing $z(s)$ at each time $s$ and player I chooses a corresponding strategy $\alpha$, then
$$
\dot x(s)=f(x(s),\alpha(z)(s), z(s)) \quad \text{for a.e $t\in (0,t)$}.
$$
In both situations,  player I wants to minimize the final payoff $g(x(t))$ and player II aims to maximize $g(x(t))$. Assume that both players play optimally, 
Theorem 4.1 in \cite{ES1984} says that
$$
u(x,t)=\sup_{\beta\in \Delta(t)}\inf_{y\in M(t)}\{g(x(t))\} \quad \text{if (\ref{eq:player1}) holds}
$$
or
$$
u(x,t)=\inf_{\alpha\in \Gamma(t)}\sup_{z\in N(t)}\{g(x(t))\} \quad \text{if (\ref{eq:player2}) holds.}
$$
In particular, if $H$ satisfies the Issacs condition \cite{I1965}:
$$
\max_{y\in Y}\min_{z\in Z}\{f(x,y,z)\cdot p\}=\min_{z\in Z}\max_{y\in Y}\{f(x,\eta,\mu)\cdot p\}, 
$$
then
$$
\sup_{\beta\in \Delta(t)}\inf_{y\in M(t)}\{g(x(t))\}=\inf_{\alpha\in \Gamma(t)}\sup_{z\in N(t)}\{g(x(t))\},
$$
which says that the game value is the same regardless which player moves first as long as both of them play optimally. 
\medskip

Compared with  the curvature G-equation, the strain G-equation has advantages (first order with a game formula in $\Rset^n$) and disadvantages (harder to discern useful strategies for estimates and leverage streamline properties in the flow). A reasonable first step towards the open question in subsection 5.2 via the differential game here is to look at flows where some knowledge of streamlines is available, such as cellular flows in 2D, ABC and Kolmogorov flows in 3D.

\section{Other Passive Scalar Combustion Models}
\setcounter{equation}{0}
In this section,  we briefly review two other ways to model the turbulent flame speeds  based on the following reaction-diffusion-advection (RDA) equation
\ba
& & T_t +  V(x)\cdot D T= d\, \Delta T+{1\over \tau_r}f(T), \no \\
& & T(x,0) = T_0(x), \;\; x\in \real^n, \label{rd1}
\ea
where $T$ represents the reactant temperature, $D$ is the spatial gradient operator, $V(x)$ is a prescribed fluid velocity, $d$ is the molecular diffusion constant,  $f$ is a nonlinear function and $\tau_r$ is the reaction time scale.  Among different nonlinearities, a popular choice is Fisher-Kolmogorov-Petrovsky-Piskunov (FKPP) reaction.
A prototypical example is $f(T)={T(1-T)}$, see \cite{Xin_00} for more details.  The reaction term $f$ corresponds to the $L^1$ term $s_l|Du|$ in G-equation model.  

\medskip

1. 
  In (\ref{rd1}), the effective burning velocity along any unit vector $p \in \Rset^n$ 
  is the minimal spreading speed $c_{p}^{*}$ from nonnegative initial data.
Similar to G-equation, two active research topics are (1) the existence of $c_{p}^{*}$ and (2) how it depends on the flow intensity if $V$ is scaled to $AV(x)$. We refer to \cite{FG_79,Fr_book,Xin_92,Xin_00,CKOR_00,BH_02,WS_04,nolen2005existence,Caff_06,NR2007,NolenXin_09,Xin_09,WS_10,
Z_10,BH_12,NRRZ_12,AZ_13,SXZ2dkpp,SXZ3dkpp, NaRo_15,IPM_2022,DP_22} and references therein, among others. A comparison between (\ref{rd1}) and G-equation was studied in \cite{XY_14b}.
 
2. 
An inviscid model was derived \cite{MS_94} from (\ref{rd1}) with KPP reaction in the scaling regime 
$$
\kappa=d\, \epsilon, \ \tau_r = \epsilon \quad \mathrm{and} \quad \ V=V\left({x\over \epsilon ^{\alpha}}\right)
$$ 
for $\alpha\in (0,1)$ and a small $\epsilon>0$ by assuming that the flame thickness is much smaller than the turbulence scale. Then the turbulent flame speed along direction $p$ is given by
$$
c_T(p)=\inf_{\lambda>0}{f'(0)+{\overline H} (p\lambda)\over \lambda}.
$$
Here $\overline H(p)$ is the effective Hamiltonian associated with the following 
 cell problem: for each $p\in \Rset^n$,  there is a unique number $\overline H(p)$ such that the inviscid quadratic Hamilton-Jacobi equation below has a periodic viscosity solution $w(x) \in C^{0,1}(\Bbb T^n)$ 
\be
d\, |p + Dw|^{2} + V(x)\cdot (p + Dw) = \overline H(p) \quad \text{in $\Rset^n$}.\label{vfcell}
\ee
We refer to  \cite{EMS} and  \cite{XY_13} for comparison with the G-equation model.  \medskip

If the flame is burning through a material at rest (i.e. $V=0$) and $f(r)={1\over \epsilon}g({r\over \epsilon})$ for a suitable non-negative function $g\in C_{0}^{1}(0,1)$,  the reaction-diffusion equation (\ref{rd1}) tends to an interesting free boundary problem as $\epsilon\to 0$ (\cite{CV1995})
\be\label{eq:C-model}
\begin{cases}
u_t-\Delta u=0 \quad \text{in $\{u>0\}$}\\[3mm]
|Du|=1 \quad \text{in $\partial \{u>0\}$}.
\end{cases}
\ee
Existence and regularity of the solution $u$ have been studied in \cite{CV1995} among others. Moreover,  this singular limit and the corresponding homogenization problem has been considered in \cite{Caff_04, Caff_06} for a generalized version when (\ref{eq:C-model}) has an advection term  in the homogeneous media. Also, see \cite{Caff_20} for a survey on other scaling related to the stationary version of the above equation and phase transition problems.

\section{Basic Stochastic G-equation}
\setcounter{equation}{0}
The existence of average front speeds in the basic G-equation 
\be
G_t+|DG|+V(x,\omega)\cdot DG=0, \label{stoch-geq}
\ee
with a stochastic incompressible vector field $V(x,\omega)$ (where $\omega$ is the random sampling variable) has been studied via homogenization 
\cite{NN_11,CS_13,BIN_20}.
Theoretical and numerical studies in physics and engineering literature include \cite{Y_88,Siv_88, KA_92,CY_98,KAW1988,ZR_94} among others.
\medskip

Due to lack of compactness, stochastic flows or random media have extreme behavior and cause additional  difficulty on the large scale so that homogenization can fail even in the coercive and convex Hamilton-Jacobi equations (\cite{Xin_09}, chapter 4 and references therein). Hence technical assumptions on the random field  are needed.  Stochastic homogenization of Hamilton-Jacobi equation in several space dimensions was first proved in  \cite{S1999} and \cite{RT2000} in the coercive, convex and stationary ergodic setting.   The strategy is to first find suitable subadditive quantities, then apply the sub-additive ergodic theorem instead of relying on corrector problems to obtain an effective quantity (often equal to the dual of the effective Hamiltonian), and finally recover the effective Hamiltonian via duality. In statistical mechanics, such a method was developed in the study of first passage percolation through the connectivity function to establish the existence of limiting shapes \cite{K_82,ACC_90}, which can be viewed as a discrete propagation problem associated with the Hamiltonian $H(p,x)=a(x,\omega)|p|$. This method, however, does not work for a non-convex Hamilton-Jacobi equation that requires completely different approaches and might fail in general \cite{Z2017, FS2017}. To the best of our knowledge,  the only available method to homogenize non-convex Hamilton-Jacobi equations under the general stationary ergodic setting is to first identify the shape of the effective Hamiltonian and then build customized correctors  \cite{ATY2015,ATY2016, G2016, QTY2018}.  When i.i.d type assumptions are imposed, approaches from first passage percolation can be used \cite{AC2015, FS2017}. We would like to mention that the periodic setting is a special case of the general stationary ergodic situation, but periodicity is on the opposite side of decorrelation. 

\subsection{Homogenization Analysis}
As in the case of periodic flows, the control formulation and the reachability concept play an essential role.  The key quantity is 
the minimal travel time $\tau(x,y,\omega)$ for the control trajectory $\xi_{\alpha}(s)$ to connect two points 
$x=\xi_{\alpha}(0)$ and $y=\xi_{\alpha}(\tau)$. In two dimensions, the incompressible flow is given by 
$V=(-\Psi_y,\Psi_x)$, $\Psi$ a stream function. Under stationarity, ergodicity,  and a third order moment condition of 
$\Psi$ (i.e. $\mathbb{E}[|\Psi(0)|^3] < +\infty$), the scaled $\tau(0, r\, y, \omega)$ as a function 
of $r \in \mathbb{R}$ is sub-additive, and 
$r^{-1}\, \tau(0, r\, y, \omega)$
converges at large $r$ to a deterministic 
function $\overline{q}=\overline{q}(y)$ almost surely by the sub-additive ergodic theorem \cite{NN_11}. Then homogenization similar to Theorem \ref{mainthm} holds with probability one based on the control formula (\ref{controlform}). The effective Hamiltonian is recovered from $\bar q$ by duality 
(Legendre transform):
\be
\bar H(p) = \sup \{ p \cdot y: y \in \mathbb{R}^d, \bar{q}(y) = 1 \}. \label{duality_barH}
\ee
That $\bar H$ being convex and homogeneous of degree one 
follows from (\ref{duality_barH}), as it is the supremum of a family of linear functions of $p$. 
\medskip

In dimensions higher than two, results in \cite{NN_11} extend under certain assumptions of $\tau(x,y,\omega)$ when $|x- y|$ diverges.
Through a refined analysis of $\tau$ and its long time behavior along a selected random sequence of times, \cite{CS_13} removed these assumptions and proved homogenization under $V=V(x,\omega)$ being stationary, ergodic, 
divergence free and zero mean. 
The sub-additive techniques extend to problems involving   
time-dependent flows.
For space-time flows $V=V(t,x,\omega)$, with a quantitative estimate of reachability of any two points by a controlled path,  
\cite{BIN_20} proved the almost sure homogenization of (\ref{stoch-geq}) under uniform bound ($\sup_{t,x,\omega}\, V < \infty$), stationarity, finite time decorrelation, small mean, Lipschitz continuity and incompressibility.  
As in periodic flows of section 3, the reachability under control in the stochastic flows  is both ways for any two points. The relatively strict finite time decorrelation assumption in \cite{BIN_20} has been relaxed to allow the infinite temporal range of dependence of the space-time flow $V$ in a recent work \cite{ZZ2023} based on a non-autonomous
sub-additive theorem.
\medskip

The travel time approach in \cite{NN_11,CS_13} is fully Lagrangian and bypasses the corrector problem entirely. It relies on  duality to come back to $\bar H$. 
The duality is lost however in a  non-convex Hamilton-Jacobi equation (e.g. when the strain rate is considered in G-equation), and it is open what would replace duality for a fully Lagrangian (corrector free) method to work.  

\begin{ques}
Can homogenization be established for strain G-equation for some interesting stochastic flows?  A basic example is 3D random  shear flow.
\end{ques}


\subsection{Formal Analysis, Numerical and Physical Experiments}
\begin{figure}
    \centering
    \includegraphics[scale=0.3]{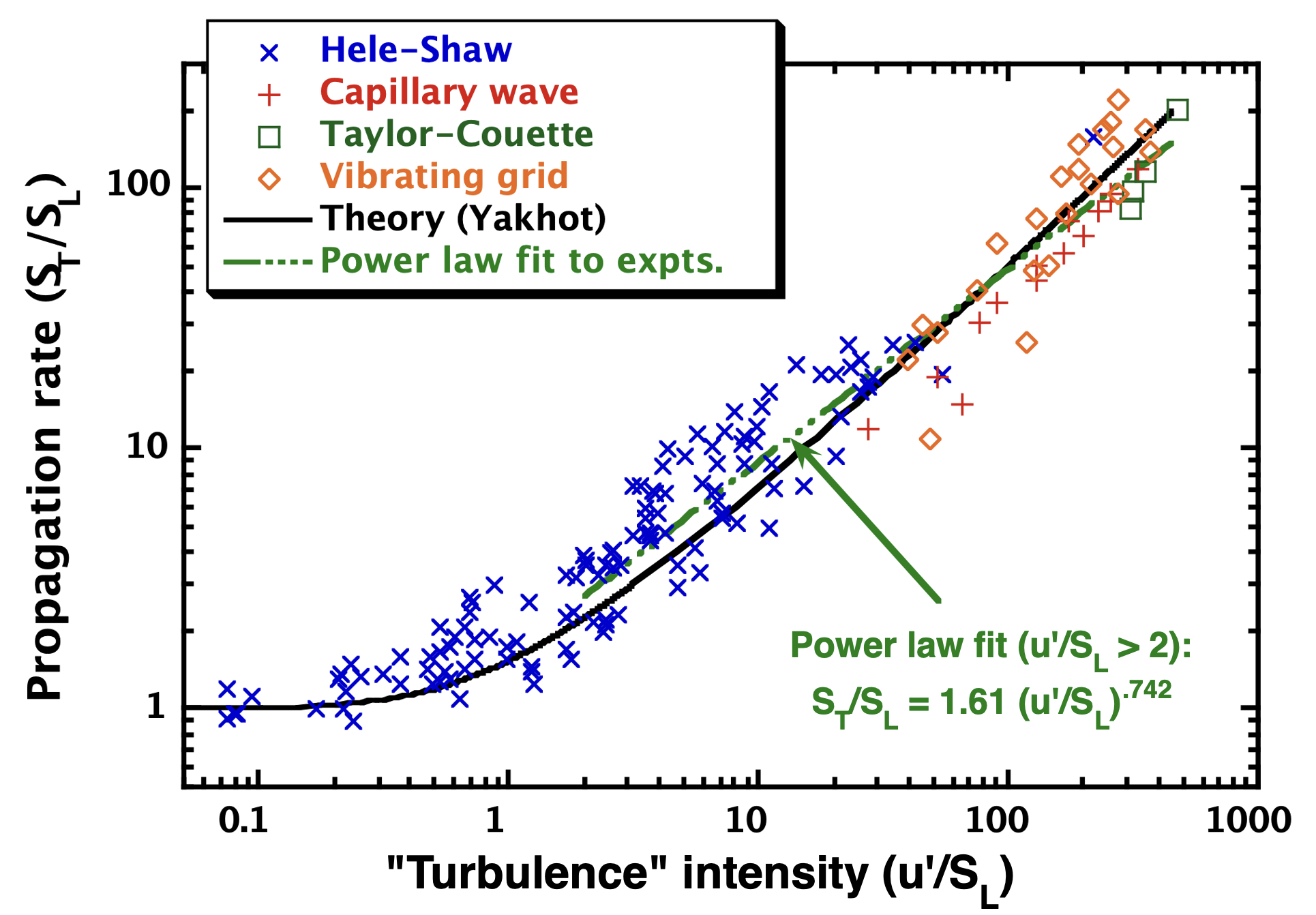}
    \caption{Comparison of propagation speeds ($s_T/s_l$) of wrinkled aqueous autocatalytic reaction fronts in Hele-Shaw, capillary wave, Taylor-Couette and vibrating-grid flows \cite{Abid_99,Ronney_95,Shy_96} along with a prediction by \cite{Y_88} and a sublinear power-law fit (applied only to data at $u'/s_l > 2$) in the range $[0.1,1000]$ of scaled flow intensity.} \label{ronney}
\end{figure}

Formal renormalization group and scaling analysis \cite{Y_88,CY_98} for the basic stochastic G-equation (\ref{stoch-geq}) showed 
that $\bar{H}=\bar{H}(A)$ scales as $O(A/\sqrt{\log A})$ at large $A$ 
for a random incompressible flow of the form  $A\, V(x,t,\omega)$ in three dimensions. 
In aqueous autocatalytic chemical reaction experiments on four different multi-scale flows, namely Hele-Shaw, capillary wave, Taylor-Couette, and vibrating-grid flows \cite{Ronney_95,Shy_96,Abid_99}, 
the theory \cite{Y_88} and a sub-linear power law fit reach good agreement with measurements in the range of flow intensities  scaled by $s_l$ (denoted by $u'/s_l$, $u':=A$ is the flow intensity) in Fig. \ref{ronney}.
\medskip

In view of the asymptotic law 
$O(A/\log A)$ in BC (cellular) flows \cite{XY_14b} and 
the law $O(A)$ in ABC and Kolmogorov flows \cite{XYZ,KLX_21}, we see that the presence of continuously many scales in stochastic flow \cite{Y_88,CY_98} indicates the absence of structures such as closed streamlines (in BC/cellular flows) and ballistic orbits (in ABC and Kolmogorov flows) as a plausible explanation of the $O(A/\sqrt{\log A})$ enhancement. A rigorous study awaits to be carried out to discover the mathematics of this long-standing random phenomenon. 
\medskip

Curvature dependence of flame speeds remain an active research topic in combustion science 
today \cite{EnsKin_23}. With the increasing concern of global warming, 
new combustion systems aim specifically
to reduce or eliminate greenhouse gas emissions. Since carbon dioxide is a key greenhouse factor, fuels without carbon content, such as hydrogen and ammonia, are promising alternatives yet with new challenges such as widened range of flammability, high flame speed and low activation energy for ignition \cite{Hydroflame_09}. 
In \cite{EnsKin_23}, premixed flames subject to turbulent disturbances and harmonic oscillation are studied experimentally and computationally. These flames come from mixtures of methane/hydrogen/air
and are stabilized on an oscillating 
holder to consistently produce harmonic wrinkles for studying curvature effects. The direct numerical simulations are also conducted for a range of 
flow intensities. The turbulent flame speed and the Markstein length are characterized via the basic G-equation as a reduced order model. A power law relating turbulent flame speed to the ensemble-averaged flame curvature in the curved flame regions is observed in \cite{EnsKin_23}. 
\medskip

Likewise,  flow induced flame stretch, a fundamental quantity measuring the rate of change of flame surface area, continues to receive attention from experimental measurements to direct simulations, see \cite{Stretch_21} and references therein.

\section{Conclusions}
We have reviewed mathematical tools for averaging (homogenizing) 
level-set G-equations arising in turbulent combustion. As curvature and flame stretch effects are included, the G-equations become non-convex and non-coercive. We showed how to leverage Lagrangian ideas (e.g. reachability of control and game trajectories) to compensate for the lack of compactness in the Eulerian (viscosity solution) framework to prove the existence and analyze qualitative properties of average flame speeds in prototypical (shear, BC, ABC and Kolmogorov) flows. We discussed  
growth laws of flame speeds as flow intensities 
become large based on the streamline structures of the ambient fluid flows. We explored physical background, and 
proposed open problems for further research in this emerging area of game theoretic analysis of non-convex and non-coercive geometric PDEs in complex advection.

\section{Acknowledgements}
This work was partially supported by NSF grants DMS-1952644, DMS-2000191, DMS-2309520, and AFOSR grant FA9550-16-1-0197. The authors would like to thank 
Prof. Y. Liu at the National Cheng-Kung University at Tainan, Taiwan, for his assistance with numerical computations.




\bigskip

 \bibliographystyle{plain}
 \bibliography{references.bib}





\end{document}